\documentclass[pdflatex,sn-mathphys-num]{sn-jnl}

\usepackage{graphicx}%
\usepackage{multirow}%
\usepackage{amsmath,amssymb,amsfonts}%
\usepackage{amsthm}%
\usepackage{mathrsfs}%
\usepackage[title]{appendix}%
\usepackage{xcolor}%
\usepackage{textcomp}%
\usepackage{manyfoot}%
\usepackage{booktabs}%
\usepackage{algorithm}%
\usepackage{algorithmicx}%
\usepackage{algpseudocode}%
\usepackage{listings}%

\theoremstyle{thmstyleone}%
\newtheorem{theorem}{Theorem}[section]
\newtheorem{lemma}[theorem]{Lemma}
\newtheorem{corollary}[theorem]{Corollary}

\theoremstyle{thmstyletwo}%
\newtheorem{remark}{Remark}%

\theoremstyle{thmstylethree}%
\newtheorem{definition}{Definition}%

\newcommand{\vp}{\varphi}
\newcommand{\HH}{\mathcal{H}}
\newcommand{\R}{\mathbb{R}}
\newcommand{\SSS}{\mathbb{S}}
\newcommand{\N}{\mathbb{N}}
\newcommand{\Z}{\mathbb{Z}}

\newcommand{\loc}{\mathrm{loc}}
\DeclareMathOperator{\sgn}{sgn}
\DeclareMathOperator{\spt}{spt}
\DeclareMathOperator{\sing}{sing}
\DeclareMathOperator{\reg}{reg}
\DeclareMathOperator{\graph}{graph}
\DeclareMathOperator{\id}{id}
\begin{document}

\title[]{On polynomial solutions to the minimal surface equation}


\author*[1]{\fnm{Yifan} \sur{Guo}}\email{yifag15@uci.edu}



\affil*[1]{\orgdiv{Department of Mathematics}, \orgname{University of California, Irvine}, \orgaddress{\street{Rowland Hall, UCI}, \city{Irvine}, \postcode{92697}, \state{CA}, \country{USA}}}



\abstract{We are interested in finding a nonlinear polynomial $P$ on $\R^n$ that solves the minimal surface equation. Even though no explicit solution is found in this article, we investigate constraints that a polynomial solution must obey.
	
	We first prove a structure theorem on such polynomials. We show that the highest degree term $P_m$ must factor as $p^kQ_m$ where $k$ is odd, $p $ is irreducible, and $Q_m\ge 0$ on $\R^n$ with $\{Q_m=0\}\subset\{p=0\}\cap\{\nabla p=0\}$. Moreover, the level sets of $P_m$ are all area-minimizing and the unique tangent cone of $\graph P$ at infinity is $\{p=0\}\times\R$. If $k\ge 3$, we know further that lower order terms down to some degree are divisible by $p$.
	
	We also show that $P$ must contain terms of both high and low degree. In particular, it cannot be homogeneous. 
	
	As a consequence of the structure theorem, we get degree estimates for polynomial solutions. We have $\deg P\ge 4$ by ruling out cubic polynomial solutions.
	Using an extended eigenvalue estimate on the Jacobi operator by Zhu \cite{zhu2018first}, we are able to show that $\mu_n^-< \deg p +k^{-1}\deg Q_m< \mu_n^+$ where $\mu_n^\pm=\frac{n-1\pm\sqrt{(n-3)^2-4(n-2)}}{2}$.
	
	Finally, we prove that $\{p=0\}$ cannot be an isoparametric minimal cone. We also show that for a nonlinear polynomial solution on $\R^8$, we have $\deg p=3$ and that $\{p=0\}$ is an area-minimizing but not strictly minimizing cone in $\R^8$. These results give strong restrictions on possible polynomial solutions to the minimal surface equation.}




\maketitle

\section{Introduction}

\subsection{Entire solutions to the minimal surface equation}
Let $P$ be a function on $\R^n$. The graph of $P$ having mean curvature $0$ in $\R^{n+1}$ is equivalent to $P$ satisfying the minimal surface equation:
\begin{align}
	\label{mse}
	(1+|\nabla P|^2)\Delta P-\sum_{i,j=1}^n  P_{x_i}  P_{x_j}  P_{x_ix_j}=0.
\end{align}
Such a function $P$, defined on all of $\R^n$, is called an entire solution to the minimal surface equation.

A well-known theorem of Bernstein \cite{B} says that the only entire solutions to the minimal surface equation on $\R^2$ are affine functions. The theorem was generalized to $\R^7$ by the works of Fleming \cite{fleming1962oriented}, De Giorgi \cite{de1965estensione},  Almgren \cite{almgren1966some} and Simons \cite{simons1968minimal}.

On the other hand, Bombieri, De Giorgi, and Giusti \cite{bombieri1969minimal} constructed nonlinear entire solutions to the minimal surface equation on $\R^n$ for $n\ge 8$ modeled on Simons' cone. Thus the restriction on $n\le 7$ in Bernstein theorem is sharp. Many other entire solutions have been constructed by Simon \cite{simon1989entire} modeled on isoparametric minimal cones (cf. Subsection \ref{isopara}).

We are interested in finding constraints on entire solutions to the minimal surface equation (\ref{mse}) which are \textit{polynomials}, a question originally raised in \cite{simon1997minimal}. There are several motivations for studying this problem.

Firstly, even though nonlinear entire solutions are known to exist, none of them are explicit. In order to find an explicit example, we should try some special ansatz and polynomials are simple and natural choices.

Secondly, polynomials fit the expected growth condition of an entire solution. Indeed, a long standing conjecture asks whether all entire solutions to the minimal surface equation have polynomial growth (see, for example \cite{bombieri1972harnack,simon1997minimal,Yau}). This conjecture is verified by Simon \cite{simon1989entire} under some technical assumptions.

Finally, in the case of the minimal surface system, where the function is vector-valued, the same problem is trivial: by a calibration argument \cite[5.4.19]{Fed}, graphs of holomorphic functions, in particular complex polynomials, are always area-minimizing and hence they provide a huge number of polynomial solutions to the system. However, in the codimension one case, the problem becomes very subtle in that there is no known example of a real polynomial whose graph is minimal except affine functions.

The analysis developed here reveals unexpected constraints arising from the interplay between nonlinear PDEs, algebraic geometry, and geometric measure theory.

\subsection{Elementary properties of polynomial solutions}
We start by fixing the terminology we will be using throughout the article.

\begin{definition}\label{polynomial.solution}
	A polynomial $P$ on $\R^n$ is called a \textit{polynomial solution} if it solves the minimal surface equation (\ref{mse}) and $P$ has degree $\ge 2$. Adding a constant to $P$ if necessary, we will always assume that $P$ has no constant term. 
\end{definition}
Then we can write a polynomial solution as the sum of homogeneous polynomials 
\begin{align}
	\label{graded}
	P=P_m+P_{m-1}+\cdots+P_2+P_1
\end{align} where $P_i$ is the degree $i$ homogeneous component of $P$. By definition, $P_m\neq0$ is the highest degree term of $P$ and $\deg P=m\ge 2$.

We observe that the left hand side of the minimal surface equation (\ref{mse}) is a polynomial and thus we can consider terms in each degree separately. To proceed, let $L$ denote the operator: \[LP:=|\nabla P|^2\Delta P-\sum_{i,j=1}^nP_{x_i}P_{x_j}P_{x_ix_j}.\] 
Then the minimal surface equation becomes
\[
\Delta P+LP=0.
\]
If the degree of $P$ is $m$, then the degree of $LP$ is $3m-4$ and the degree of $\Delta P$ is $m-2$. Thus equation (\ref{mse}) is equivalent to $3m-3$ equations for $P_m,\dots,P_1$ which is obtained by expanding (\ref{mse}) according to (\ref{graded}). The equation of the highest degree $3m-4$ is
\begin{equation}
	LP_m=0.\label{1lap0}
\end{equation}
Note that (\ref{1lap0}) is equivalent to the fact that every level set of $P_m$ has mean curvature 0 at each regular point since the mean curvature of the level set at such a point is given by $|\nabla P_m|^{-3}L P_m$. For more details, see Lemma \ref{expansion}.

Solving $LP_m=0$ is the first step in solving the full minimal surface equation and the only known solutions which are polynomials are powers of linear functions. That powers of linear functions cannot be the leading term of a polynomial solution (of degree $\ge 2$) can be seen from Theorem \ref{struct0} (ii).
One nontrivial solution to $Lu=0$ is the function $F$ constructed in \cite[Section III]{bombieri1969minimal} which leads to the first proof of area-minimality of Simons' cone. However, by Theorem \ref{noisopara0}, this function as well as its powers cannot be the leading term of a polynomial solution.

Using equation (\ref{1lap0}), we have the following elementary result showing that polynomial solutions must contain terms of both high and low degree.

\begin{theorem}(Theorem \ref{homo})\label{homo0}
	Let $P=P_m+\cdots+P_1$ be a polynomial solution. Then there exists $1\le i\le \lceil\frac{m}{3}\rceil$ such that $P_i\neq0$. Here $\lceil x\rceil$ is the smallest integer greater than or equal to $x$. In particular, there are no homogeneous polynomial solutions.
\end{theorem}

The same method will give a slight variant of the theorem which rules out quadratic polynomial solutions.
\begin{theorem}(Theorem \ref{nodeg2})\label{nodeg20}
	Polynomial solutions cannot be of the form $P=P_m+P_1$. In particular, there are no quadratic polynomial solutions.
\end{theorem}

\subsection{The main structure theorem of polynomial solutions}
We now state our main structure theorem of polynomial solutions.
The notation and concepts in geometric measure theory plays an important role in both the statement and the proof of the theorem. We will recall the necessary background in Section \ref{notations}. 
\begin{theorem}(Theorem \ref{struct})\label{struct0}
	Let $P=P_m+\cdots+P_1$ be a polynomial solution. Then we have:
	\begin{enumerate}
		\item[(i)] For every $t\in\R$, the current $\partial [P_m<t]$ is area-minimizing. If $t\neq0$, $\partial[P_m<t]$ is a smooth minimal hypersurface.
		\item[(ii)] There exist 
		\begin{enumerate}
			
			\item [(a)]$k\ge 1$ odd, 
			
			\item[(b)]$p $ an irreducible homogeneous polynomial that changes sign and has $\deg p\ge3$,
			
			\item[(c)] $Q_m$ an even degree homogeneous polynomial which is non-negative on $\R^n$ and coprime with $p $, with $\{Q_m=0\}\subset\{p=0\}\cap\{\nabla p=0\}$ of dimension $\le n-8$,
		\end{enumerate}
		so that
		\begin{align}\label{pkq0}
			P_m=p ^kQ_m.
		\end{align}
		
		\item[(iii)] The unique tangent cone of $\graph P$ at infinity is $\partial [p <0]\times \R$.
		
		\item[(iv)]The function $u=\frac{1}{|\nabla (p Q_m^{1/k})|}$ is the velocity of the family of minimal hypersurfaces $(\{P_m=t^k\})_{t\in\R}$ and is finite on the regular part of $\{P_m=t^k\}$ for each $t\in\R$.

		\item[(v)]  If $k\ge 3$ in (\ref{pkq0}), then there exists $2\le s\le m-1$ so that
		\[P=p ^kQ_m+p ^{k_{m-1}} Q_{m-1}+\cdots +p ^{k_{s+1}}Q_{s+1}+p Q_s+P_{s-1}+\cdots+P_1\]
		where  $p $ and $Q_i$ are coprime for $s\le i\le m$ and $k_i\ge 2$ for $s+1\le i\le m-1$.
	\end{enumerate}
\end{theorem}
We remark here that (i) aligns well with the fact that $P_m$ solves (\ref{1lap0}). In view of (ii) and (iii) we make the following definition.
\begin{definition}
	Let $P$ be a polynomial solution with $P_m=p^kQ_m$ as in (\ref{pkq0}). Then the polynomial $p$ is called the \textit{tangent cone factor} of $P$.
\end{definition}

The main technique used to prove Theorem \ref{struct0} is to consider the limit of translated blow-down: we scale down $\graph P$ by $\lambda$ and translate vertically by $t\lambda^{-m+1}$ for some fixed $t\in\R$ (see also \cite[Proof of Theorem 4]{simon1989entire}). We are able to show the following convergence of currents in Lemma \ref{convgr} and Theorem \ref{struct} (iii) \[\lim_{\lambda\to 0}\eta_{\lambda\#}\graph P-t\lambda^{-m+1}e_{n+1}= \partial[P_m<t]\times\R\]
where $\eta_\lambda:\R^n\to \R^n$ is the map $x\mapsto \lambda x$.
The structure theorem then follows from the properties of area-minimizing currents and the real Nullstellensatz.

\subsection{Degree estimates of polynomial solutions}

Using Theorem \ref{struct0}, we can get some estimates for the degree of a polynomial solution. First, as a direct consequence, we can rule out cubic polynomial solutions.
\begin{theorem}(Theorem \ref{nodeg23})\label{nodeg30}
	Let $P$ be a polynomial solution. Then $\deg P\ge 4$.
\end{theorem}

Next, we can estimate the degree of the tangent cone factor $p$.

\begin{theorem}(Theorem \ref{n-2})\label{n-20}
	Let $P$ be a polynomial solution  on $\R^n$ and $P_m=p^kQ_m$ as in Theorem \ref{struct0} (ii). Then 
	\[\mu_n^-< \deg p +k^{-1}\deg Q_m< \mu_n^+\] where
	$\mu_n^{\pm}= \frac{n-1\pm\sqrt{(n-3)^2-4(n-2)}}{2}$. In other words, $k\mu_n^-<\deg P< k\mu_n^+$.
\end{theorem}
It would be interesting to get rid of the dependence on $k$ in the estimate for $\deg P$. However, it seems that this is the best we can do with our approach.

Our proof utilizes Theorem \ref{struct0} (iv) as well as an extension of Simons' eigenvalue estimate to singular setting by Zhu \cite{zhu2018first}.

\subsection{Further restrictions on the tangent cone of \texorpdfstring{$\graph P$}{graph P}}
Theorem \ref{n-20} shows that the degree of the tangent cone factor $p$ is bounded by a constant depending only on the dimension $n$. In this section, we explore further restrictions on the tangent cone factor $p$ imposed by Theorem \ref{struct0}. These constraints are strong enough to exclude all currently known examples of area-minimizing cones from occurring as the tangent cone of $\graph P$.

\label{isopara}

To state the first result, we recall that a minimal hypercone $C\subset\R^n$ is called \textit{isoparametric} if the link of the cone $C\cap \SSS^{n-1}$ has constant principal curvatures in $\SSS^{n-1}$ \cite{Car0}. 
To the author's best knowledge, isoparametric minimal cones and their products with $\R^l$ are the only known examples of area-minimizing hypercones. Note that these include Lawson's cones, which are cones over products of two spheres. Every isoparametric minimal cone is the zero locus of a homogeneous polynomial whose degree equals that of the Cartan–Münzner polynomial associated with the isoparametric family \cite{Car,Mu1,Mu2}. This fact is a direct consequence of Lemma~\ref{tilde.p}. The first restriction we establish rules out the possibility of isoparametric minimal cones as tangent cones for polynomial solutions.
\begin{theorem}
    (Corollary \ref{noisopara.geometric})\label{noisopara0}
	Let $P$ be a polynomial solution. Then the tangent cone at infinity of $\graph P$ cannot be $(C\times\R^{l-1})\times \R$ where $C$ is an isoparametric minimal cone and $l\ge 1$.
\end{theorem}

To motivate the second result, we recall the Bernstein's theorem which says that the smallest dimension in which a nonlinear entire solution could exist is $n=8$. For polynomial solutions on $\R^8$, we have the following theorem on its tangent cone.

\begin{theorem}(Theorem \ref{R8cubic})\label{R8cubic0}
	Let $P=P_m+\cdots+P_1$ be a polynomial solution on $\R^8$ with $P_m=p^kQ_m$ as in Theorem \ref{struct0} (ii). Then $\deg p=3$ and $\{p=0\}\subset \R^8$ is a strictly stable, area-minimizing cone which is not strictly minimizing (for a definition, see Definition \ref{sssm}).
\end{theorem}
We can rephrase it as a Bernstein type theorem.
\begin{corollary}\label{R8}
	If $P$ is a polynomial on $\R^8$ solving the minimal surface equation and the tangent cone of $\graph P$ at infinity is $C\times\R$ where $C$ is strictly minimizing, then $P$ is an affine function.
\end{corollary}
Theorem \ref{R8cubic0} uses the growth estimates of Simon \cite{simon1989entire} and the fact that no cubic solutions exist from Theorem \ref{nodeg30}. The conclusion that $\{p=0\}$ is not strictly minimizing comes from the assumption of \cite{simon1989entire}. Such an assumption appears in many results of minimal hypersurfaces. Note that there are no known examples of an area-minimizing but not strictly minimizing hypercones except for $\R^2\subset \R^3$ (see \cite{lin_1987}). Moreover, in $\R^8$, the only known area-minimizing cones are cones over $S^3\times S^3$ and $S^2\times S^4$, both of which are isoparametric and strictly minimizing. Finally, we note that Tkachev \cite{TkClifford} constructed a family of cubic minimal cones via Clifford algebras in various dimensions, and in particular in $\R^8$. He also classified the so called radial cubic minimal cones in \cite{tkachev2010classification}, see also the book \cite[Chapter 6]{NTV}. However, the examples in \cite{TkClifford} in $\R^8$ are not area-minimizing since they have 2-dimensional singular sets.

The outline of the rest of the paper is as follows. We first introduce some notation and background in geometric measure theory and real algebraic geometry in Section \ref{notations}. In Section \ref{2}, we discuss some elementary properties of polynomial solutions. In Section \ref{3}, we prove the main structure theorem of polynomial solutions.  In Section \ref{4}, we use the structure theorem to estimate the degree of a polynomial solution. In Section \ref{5}, we prove Theorem \ref{noisopara0} and \ref{R8cubic0}.

\section{Notation and Preliminaries}
\label{notations}
\subsection{Notation}
We introduce the following notation:

$\R$: the set of real numbers;

$\N$: the set of positive integers;

$\R[x_1,\dots,x_n]$: the polynomial ring of $n$ variables over $\R$;

$p\mid q$: the polynomial $p$ divides the polynomial $q$;

$p\nmid q$: the polynomial $p$ does not divide the polynomial $q$;

$\deg p$: the degree of the polynomial $p$;

$Z(p)$: the set $\{x\in\R^n:p(x)=0\}$ for $p\in \R[x_1,\dots,x_n]$;

$(p)$: the ideal generated by $p$ in $\R[x_1,\dots,x_n]$;

$I(S)$: the ideal generated by polynomials in $\R[x_1,\dots,x_n]$ that vanish on $S\subset \R^n$; 


$|x|$: the length of the vector $x\in \R^n$;

$B_r(x)$: the set $\{y\in \R^n:|x-y|<r\}$ (write $B_r$ if $x=0$);

$\SSS^{n-1}$: the set $\{x\in \R^n:|x|=1\}$;

$d(x,K)$: the distance from a point $x\in \R^n$ to a closed set $K\subset \R^n$;

$\eta_\lambda$: the map $\R^n\to \R^n$ given by $\eta_\lambda(x)=\lambda x$ for $\lambda>0$.

For $u:\R^n\to\R$ a $C^2$ function and $\Sigma\subset \R^n$ an orientable smooth hypersurface, we introduce:

$u_{x_{i_1}\dots x_{i_k}}:=\frac{\partial^k u}{\partial x_{i_1}\cdots\partial x_{i_k}}$;

$\graph u$: the $C^2$ manifold $\{(x,u(x)):x\in \R^n\}$ in $\R^{n+1}$;

$\nabla u$: the Euclidean gradient of $u$ at $x$;

$\nabla_\Sigma u$: the projection of $\nabla u(x)$ onto $T_x\Sigma$;

$\Delta u$: Euclidean Laplace of $u$;

$\Delta_\Sigma u$: the Laplace-Beltrami operator of $\Sigma$ acting on $u$;

$|A_\Sigma|^2$: the norm squared of the second fundamental form of $\Sigma$;

$\nu_\Sigma(x)$: the unit normal of $\Sigma$ at $x$;

$\graph_\Sigma u$: the $C^2$ manifold $\{x+u(x)\nu_\Sigma(x):x\in \Sigma\}$ for $\|u\|_{L^\infty}$ sufficiently small.

The following are fundamental concepts in geometric measure theory. We refer to \cite{Fed,morgan2016gmt,Sim} for their precise definitions.

For $U$ an open set in $\R^n$, we introduce:

$R_k^\loc(U)$: the set of integer multiplicity $k$-rectifiable currents in $U$;

$S_j\rightharpoonup S$: weak convergence of the currents $S_j\in R_k^\loc(U)$ to $S\in R_k^\loc(U)$ as $j\to\infty$;

$\|S\|$: the Radon measure associated to $S\in R_k^\loc(U)$;

$S\llcorner V$, $\mu\llcorner V$: restriction of the current $S\in R_k^\loc(U)$ or the Radon measure $\mu$ to a set $V\subset U$;

$\partial S$: the boundary of $S\in R_k^\loc(U)$;

$\spt S$: the support of $S\in R_k^\loc(U)$;

$\reg S$: the regular set of $S\in R_k^\loc(U)$;

$\sing S$: the singular set of $S\in R_k^\loc(U)$;

$I_k^\loc(U)$: the set of $S\in R_k^\loc(U)$ such that $\partial S\in R_{k-1}^\loc(U)$;

$[M]$: the $k$-current associated to a $k$-dimensional smooth submanifold $M\subset \R^n$;

$[f<0]$: the $n$-current defined by the open set $\{f<0\}\subset\R^n$ where $f:\R^n\to\R$ is continuous;

$\phi_\#S$: the push forward of a current $S\in R_k^\loc(U)$ via a proper smooth map $\phi:U\to V$ where $V\subset \R^m$ is open.

\subsection{Area-minimizing boundaries}
Since area-minimizing currents play an important role in the main structure theorem, Theorem \ref{struct}, we recall its definition and relevant properties.
\begin{definition}[\cite{Fed,morgan2016gmt,Sim}]
	Let $U\subset \R^n$ be open. A current $T\in R_k^\loc(U)$ is \textit{area-minimizing} in $U$ if for any compact subset $K\subset U$ and any $S\in R_k^\loc(U)$ satisfying $\partial S=\partial T$ and $\spt(T-S)\subset K$, we have $\|T\|(K)\le \|S\|(K)$.
	
	If $T\in I_{n-1}^\loc(\R^n)$ is area-minimizing and $T=\partial [V]$ for some open set $V\subset \R^n$, then $T$ is called an \textit{area-minimizing boundary}.
\end{definition}

The area-minimizing property comes naturally with entire solutions.
\begin{lemma}[{\cite[\S 3]{simon1983survey}}]\label{minimizing}
	Let $u:\R^n\to\R$ be an entire solution to the minimal surface equation. Then $[\graph u]\in I_n^\loc(\R^{n+1})$ is an area-minimizing boundary in $\R^{n+1}$.
\end{lemma}

The following two lemmas are the properties of area-minimizing boundaries that we need.
\begin{lemma}[{\cite[Theorem 37.2]{Sim}}]\label{FedFlem}
	Let $U\subset \R^n$ be open. Suppose $\{T_j\}\subset I_{n-1}^\loc(U)$ is a sequence of area-minimizing boundaries in $U$. Then there exist a subsequence $\{j'\}\subset \{j\}$ and $T\in I_{n-1}^\loc(U)$ such that $T_{j'}\rightharpoonup T$ in $I_{n-1}^\loc(U)$. Moreover, $T$ is area-minimizing in $U$ and $\|T_{j'}\|\to \|T\|$ as Radon measures in $U$.
\end{lemma}

\begin{lemma}[{\cite[Theorem 1]{bombieri1972harnack}}]\label{indecomposable}
	Let $\partial [U]$ be an area-minimizing boundary in $\R^{n+1}$. Then $U$ is connected.
\end{lemma}

We also recall here the notion of tangent cone at infinity of an area-minimizing boundary in $\R^n$.
\begin{definition}
	Let $T\in I_{n-1}^\loc(\R^n)$ be an area-minimizing boundary. We say that $C\in I_{n-1}^\loc(\R^n)$ is a \textit{tangent cone of $T$ at infinity} if there is a sequence $\lambda_j\to0$ such that $\eta_{\lambda_{j}\#}T\rightharpoonup C$ in $I_{n-1}^\loc(\R^n)$. We note that $T$ always has a tangent cone at infinity by Lemma \ref{FedFlem}.
\end{definition}
\subsection{Real algebraic geometry}
Various sets will arise in the later sections, and all of them are defined by polynomial equalities and inequalities. Such sets are called semi-algebraic and admit a well-defined notion of dimension, which we now introduce. See \cite{bochnak2013real} for more details.
\begin{definition}
	A set $V\subset\R^n$ is called \textit{algebraic} if there is a polynomial $f\in \R[x_1,\dots,x_n]$ such that $V=Z(f)$.
	
	A set $V\subset\R^n$ is called \textit{semi-algebraic} if it is of the form
	\[V=\bigcup_{i=1}^s\bigcap_{j=1}^{r_i}\{x\in \R^n:f_{i,j}(x)*_{i,j}0\}\]
	where $s,r_i\in \N$, $f_{i,j}\in \R[x_1,\dots,x_n]$ and $*_{i,j}$ is either $<$ or $=$ for $1\le i\le s$, $1\le j\le r_i$.
\end{definition}
\begin{remark}
	In the definition of algebraic sets, we only need one polynomial $f$. This is because if $V=Z(f_1)\cap \cdots\cap Z(f_m)$ then $V=Z(f_1^2+\dots+f_m^2)$. This is special to the field $\R$.
\end{remark}
\begin{definition}[{\cite[page 5]{hartshorne1977algebraic}}]\label{dimension}
	Let $V$ be an algebraic set. The \textit{dimension} of $V$, denoted by $\dim V$, is the maximal length $d$ of chains $V_0\subsetneq\cdots\subsetneq V_d$ of nonempty distinct irreducible Zariski closed subsets of $V$.
	
	If $V$ is semi-algebraic, then the dimension of $V$ is the dimension of the Zariski closure of $V$ which is an algebraic set.
\end{definition}
\begin{remark}
	The definition of Zariski topology as well as irreducibility can be found in \cite[Theorem 2.8.3]{bochnak2013real} and \cite[page 2,3]{hartshorne1977algebraic}.
\end{remark}
\begin{remark}
	The dimension of $V$ can also be defined as the Krull dimension of the ring $\R[x_1,\dots,x_n]/I(V)$ (see, e.g. \cite[Definition 2.8.1]{bochnak2013real}), or the Hausdorff dimension (see e.g. \cite{morgan2016gmt}). They are all equivalent in this setting (see e.g. \cite[Proposition 1.7]{hartshorne1977algebraic} and \cite[Corollary 2.8.9, Proposition 2.9.10]{bochnak2013real}).
\end{remark}
%

For later reference, we quote here the Nullstellensatz over $\R$. 

\begin{lemma}[Real Nullstellensatz]\label{null}
	Let $p\in \R[x_1,\dots,x_n]$ be an irreducible polynomial. The following properties are equivalent.
	\begin{enumerate}
		\item[(i)]   There exists $x\in Z(p)$ so that $\nabla p(x)\neq 0$. 
		
		\item[(ii)]  As a function on $\R^n$, $p$ changes sign. 
		
		\item[(iii)]  We have $(p)=I(Z(p))$.
		
		\item[(iv)]  We have $\dim Z(p)=n-1$.
	\end{enumerate}
	
	In particular, if $p$ satisfies any of the conditions above and $q$ is a polynomial that vanishes on a semi-algebraic set $V\subset \{p=0\}$ with $\dim V=n-1$, then $p\mid q$.
\end{lemma}
\begin{proof}
	Items (i)-(iv) can be found in \cite[Theorem 4.5.1]{bochnak2013real}. We justify here the last claim. Since $\dim V=n-1$, by Definition \ref{dimension} there are $V_0\subsetneq V_1\subsetneq\cdots\subsetneq V_{n-1}\subset \bar V$ where $V_i$ are nonempty irreducible Zariski closed subsets of the Zariski closure $\bar V$ of $V$. We claim that $\{p=q=0\}=\{p=0\}$. Indeed, suppose $\{p=q=0\}\subsetneq \{p=0\}$. Since $V\subset \{p=q=0\}$, we have $\bar V\subset \{p=q=0\}\subsetneq \{p=0\}$ and hence we have a chain of algebraic subsets of $\{p=0\}$
    \begin{equation}\label{dim=n}
	    V_0\subsetneq  V_1\subsetneq\cdots\subsetneq V_{n-1}\subsetneq V_n=\{p=0\}.
	\end{equation} For $1\le i\le n-1$, since $V_i$ is an irreducible Zariski closed subset of $\bar V$ which is Zariski closed, it is also irreducible Zariski closed subsets of $\{p=0\}$. Since $p$ is irreducible and $I(Z(p))=(p)$, we have $V_n$ is irreducible and Zariski closed by \cite[Theorem 2.8.3]{bochnak2013real}. By \eqref{dim=n} and Definition \ref{dimension}, we have $\dim\{p=0\}\ge  n$ which contradicts (iv). Thus $q$ vanishes on $\{p=0\}$. By (iii), we have $p\mid q$.
\end{proof}

\begin{remark}
	If $p$ is a polynomial that does not satisfy the conditions in Lemma \ref{null}, then it cannot change sign, $\{p=0\}$ is of high codimension and $\nabla p\equiv0$ on $\{p=0\}$. One such example is the polynomial $x_1^2+x_2^2\in \R[x_1,x_2,x_3]$.
\end{remark}

\section{Elementary properties of polynomial solutions}
\label{2}
In this section, we establish some elementary properties of polynomial solutions.
\begin{lemma}\label{expansion}
	Let $P=P_m+\cdots+P_1$ be a polynomial solution. Then 
	\begin{enumerate}
		
		\item[(i)] $P_m$ satisfies $LP_m=0$ where \[Lu:=|\nabla u|^2\Delta u-\sum_{i,j=1}^nu_{x_i}u_{x_j}u_{x_ix_j}.\]
		
		\item[(ii)] For any $t\neq 0$, the set $\{P_m=t\}$ is a smooth minimal hypersurface in $\R^n$ with $\nabla P_m$ nonvanishing on $\{P_m=t\}$.
		
		\item[(iii)] For any $t\neq 0$, the function $u=\frac{1}{|\nabla P_m|}$ is the velocity of the family of hypersurfaces $(\{P_m=t\})_{t\neq0}$ and hence $u$ is a Jacobi field on $\{P_m=t\}$ for $t\neq 0$ i.e. $\Delta_{\{P_m=t\}}u+|A_{\{P_m=t\}}|^2u=0$.
		
		\item[(iv)] $P_{m-1}$ satisfies 
		\begin{equation}
			\label{m-1}\begin{aligned}
				DL(P_m)P_{m-1}:=&|\nabla P_m|^2\Delta P_{m-1}+2\Delta P_m\nabla P_m\cdot \nabla P_{m-1}\\
				-\sum_{i,j=1}^n&\left( P_{m,x_i}P_{m,x_i}P_{m-1,x_ix_j}+2P_{m,x_j}P_{m,x_ix_j}P_{m-1,x_i}\right)=0.
			\end{aligned}
		\end{equation}
		\item[(v)] For any $t\neq 0$, the function $v=-\frac{P_{m-1}}{|\nabla P_m|}$ is a Jacobi field on $\{P_m=t\}$.
	\end{enumerate}
\end{lemma}
\begin{remark}
	For item (ii), after we establish the main theorem, Theorem \ref{struct}, we can normalize the speed of the level sets $\{P_m=t\}$ so that the velocity $u$ is well defined on the regular part of $\{P_m=0\}$ (cf. Theorem \ref{struct} (iv)).
\end{remark}
\begin{remark}\label{noinfo}
	If we have $P_m$ solving $LP_m=0$, then $P_{m,x_i}$ for $1\le i\le n$ are solutions to (\ref{m-1}). This corresponds to translations of the level sets in $\R^n$. Note that (v) is a consequence of (\ref{m-1}). We are going to show in Lemma \ref{convgr} that $v$ appears as the velocity of a family of minimal hypersurfaces which are the translated blow-down of $\graph P$.
\end{remark}
\begin{proof}
	(i) Since $m\ge 2$, we have $3m-4>m-2$. The desired equation is the degree $3m-4$ term in the expansion of (\ref{mse}) according to (\ref{graded}).

	(ii) Since $P_m$ is homogeneous of degree $m$, by Euler's identity we have, for any $x\in\{P_m=t\}$ \[x\cdot\nabla P_m(x)=mP_m(x)=mt\neq 0.\] Thus $\nabla P_m(x)\neq 0$ and $\{P_m=t\}$ is a smooth hypersurface. Since the mean curvature of $\{P_m=t\}$ is given by $|\nabla P_m|^{-3}L P_m$, we have $\{P_m=t\}$ is a smooth minimal hypersurface.

	(iii) Since the family of hypersurfaces are given by level sets of $P_m$, its velocity vector field is $\frac{1}{|\nabla P_m|}$. Since $(\{P_m=t\})_{t\neq0}$ is a family of smooth minimal hypersurfaces, its velocity field $\frac{1}{|\nabla P_m|}$ satisfies the Jacobi equation on $\{P_m=t\}$.
	
	(iv) Since $m\ge 2$, we have $3m-5>m-2$. The desired equation is the degree $3m-5$ term in the expansion of (\ref{mse}) according to (\ref{graded}).
	
	(v)
	First we claim that for $t\neq 0$ fixed, $v(x)=-\frac{P_{m-1}}{|\nabla P_m|}$ is the velocity field of the family of hypersurfaces $(\{P_m+sP_{m-1}=t\})_{s\in(-\epsilon,\epsilon)}$ at $s=0$. Indeed, suppose $\gamma(s)$ is a curve in $\R^n$ so that $\gamma(0)=x\in \{P_m=t\}$, $\dot \gamma(0)=v(x)\frac{\nabla P_m(x)}{|\nabla P_m(x)|}$ and $(P_m+sP_{m-1})(\gamma(s))=t$ for all $s\in(-\epsilon,\epsilon)$. Then 
	\[\begin{aligned}
		0=&\left.\frac{d}{ds}\right|_{s=0}(P_m+sP_{m-1})(\gamma(s))=P_{m-1}(x)+\nabla P_m\cdot \left(v(x)\frac{\nabla P_m(x)}{|\nabla P_m(x)|}\right)\\
		=&P_{m-1}(x)+v(x)|\nabla P_m(x)|.
	\end{aligned}\]
	Thus $v=-\frac{P_{m-1}}{|\nabla P_m|}$ and the claim is proved.
	
	The mean curvature of the level set of $P_m+sP_{m-1}$ is \[\frac{L(P_m+sP_{m-1})}{|\nabla (P_m+sP_{m-1})|^3}.\] By (i) and (iv), we have
	\[\left.\frac{d}{ds}\right|_{s=0}\frac{L(P_m+sP_{m-1})}{|\nabla (P_m+sP_{m-1})|^3}=\frac{DL(P_m)P_{m-1}}{|\nabla P_m|^3}+L(P_m)\left.\frac{d}{ds}\right|_{s=0}\frac{1}{|\nabla (P_m+sP_{m-1})|^3}=0.\] Thus $\{P_m+sP_{m-1}=t\}$ has 0 mean curvature up to the first order in $s$. On the other hand, by the linearization of the mean curvature at a minimal hypersurface, we have
	\[0=\left.\frac{d}{ds}\right|_{s=0}\frac{L(P_m+sP_{m-1})}{|\nabla (P_m+sP_{m-1})|^3}=\Delta_{\{P_m=t\}}v+|A_{\{P_m=t\}}|^2v.\] Thus $v$ is a Jacobi field on $\{P_m=t\}$ for $t\neq 0$.
\end{proof}
\begin{theorem}\label{homo}
	Let $P=P_m+\cdots+P_1$ be a polynomial solution. Then there exists $1\le i\le \lceil\frac{m}{3}\rceil$ such that $P_i\neq0$. Here $\lceil x\rceil$ is the smallest integer greater than or equal to $x$.
	In particular, there are no homogeneous polynomial solutions.
\end{theorem}
\begin{proof} Suppose there is a polynomial solution so that $P_i=0$ for all $1\le i\le \lceil\frac{m}{3}\rceil$. Let $k$ be the smallest $k$ such that $P_k\neq0$. Then $k\ge \lceil\frac{m}{3}\rceil+1$. The lowest degree term in $LP$ is $LP_{k}=0$ which has degree $3k-4>m-2=\deg \Delta P_m$. Thus the $m-2$ degree term in $\Delta P+LP$ is $\Delta P_m$ and we have $\Delta P_m=0$. Combining with $LP_m=0$ by Lemma \ref{expansion} (i), we have a nonzero homogeneous polynomial $P_m$ satisfying the $\infty$-Laplacian equation
	\[\sum_{i,j=1}^nP_{x_i}P_{x_j}P_{x_ix_j}=0.\]
	This is a contradiction in view of the next lemma.
\end{proof}

\begin{lemma}[{\cite[Proposition 4.1]{tkachev2016non}}]\label{infhar}
	Let $P$ be a smooth  function  on $\R^{n}\setminus \{0\}$ which is homogeneous of degree $m\notin \{0,1\}$ and satisfies $\sum_{i,j=1}^nP_{x_ix_j}P_{x_i}P_{x_j}=0$. Then $P=0$.
\end{lemma}
\begin{remark}
	Note that $m$ is not required to be an integer. Hence the lemma actually shows that there is no homogeneous solution to the minimal surface equation on $\R^{n}\setminus\{0\}$.
\end{remark}
We write down the proof of the lemma for the sake of convenience.
\begin{proof}
	We can rewrite the $\infty$-Laplacian equation as 
	\begin{equation}\label{inflap}
		\frac{1}{2}\langle \nabla P,\nabla |\nabla P|^2\rangle=0.
	\end{equation} 
	Let $x\in \SSS^{n-1}$ so that $P(x)=\max _{y\in \SSS^{n-1}}P(y)$. By the method of Lagrange multipliers, we have $\nabla P(x)=\lambda x$ where $\lambda\in\R$ is some real number. Evaluating (\ref{inflap}) at $x$, and using Euler's formula for the homogeneous function $|\nabla P|^2$, we have
	\[0=\frac{1}{2}\langle\lambda x,\nabla|\nabla P|^2(x)\rangle=\lambda (m-1)|\nabla P(x)|^2=\lambda^3(m-1).\]
	Since $m\neq 1$, we have $\lambda=0$ and thus $\nabla P(x)=0$.
	
	On the other hand, using Euler's formula for the homogeneous function $P$, we have
	\[mP(x)=\langle x,\nabla P(x)\rangle=0.\]
	Hence $\max _{y\in \SSS^{n-1}}P(y)=0$. Similarly we have $\min _{y\in \SSS^{n-1}}P(y)=0$ and thus $P=0$.
\end{proof}

The following theorem is a slight variant of Theorem \ref{homo}.
\begin{theorem}\label{nodeg2} 
	Polynomial solutions cannot be of the form $P=P_m+P_1$.
	In particular, there are no quadratic polynomial solutions.
\end{theorem}

\begin{proof}
	Suppose for contradiction that $P=P_m+P_1$ is a polynomial solution. After rotating the coordinates in $\R^n$, we may assume that $P_1(x)=ax_1$ with $a\neq 0$. Plugging $P$ in (\ref{mse}), and considering degree $m-2$ and $2m-3$ terms, we obtain
	\begin{align}
		(1+a^2)\Delta P_m-a^2 P_{m, x_1x_1}=0, \label{eqi}\\
		P_{m, x_1} \Delta P_m-\sum_{i=1}^n P_{m,x_i x_1}P_{m,x_i}=0.\label{eqii}
	\end{align}
	Eliminating $\Delta P_m$ in (\ref{eqi}) and (\ref{eqii}), we get
	\[\frac{a^2}{1+a^2}P_{m, x_1x_1}P_{m, x_1} -\sum_{i=1}^n P_{m,x_i x_1}P_{m,x_i}=0\]
	i.e.
	\[\left(\frac{1}{1+a^2}P_{m,x_1}^2+P_{m,x_2}^2+\cdots +P_{m,x_n}^2\right)_{x_1}=0.\]
	Thus we have $\frac{1}{1+a^2}P_{m,x_1}^2+\sum_{i=2}^nP_{m,x_i}^2$ is independent of $x_1$. 
	
	We now arrange the terms in $P_m$ according to the degree in $x_1$:
	\[P_m(x)=c_0x_1^m +c_{1}(x^\prime)x_1^{m-1}+\cdots+c_m(x^\prime)\]
	where $x^\prime=(x_2,\dots,x_n)$ and $c_i$ is a homogeneous polynomial in $x^\prime$ of degree $i$.
	
	The term in $\frac{1}{1+a^2}P_{m,x_1}^2+\sum_{i=2}^nP_{m,x_i}^2$ with the highest $x_1$ degree is 
	\[\frac{1}{1+a^2}c_0^2m^2x_1^{2m-2}+\sum_{i=2}^n c_{1,x_i}^2x_1^{2m-2}.\]
	Thus we must have $c_0=c_1=0$.
	
	Assume that $j\ge 2$ is the smallest index such that $c_j\neq 0$, then we have $c_0=\dots=c_{j-1}=0$. The term in $\frac{1}{1+a^2}P_{m,x_1}^2+\sum_{i=2}^nP_{m,x_i}^2$ with the highest $x_1$ degree is 
	\[\sum_{i=2}^{n}c_{j,x_i}(x^\prime)^2x_1^{2m-2j}.\]
	Since $c_j\neq 0$, in order for $\frac{1}{1+a^2}P_{m,x_1}^2+\sum_{i=2}^nP_{m,x_i}^2$ to be independent of $x_1$, we must have $m=j$. Hence $c_0=\dots=c_{m-1}=0$ and $P_m$ is independent of $x_1$.
	
	By equation (\ref{eqi}), we have $\Delta P_m=0$. Since $LP_m=0$, Lemma \ref{infhar} implies that $P_m=0$ i.e. $P=P_1$. This contradicts Definition \ref{polynomial.solution}.
\end{proof}
\section{Structure of polynomial solutions}
\label{3}
\subsection{Translated blow-down}
We will be studying the translated blow-down of the graph of a general polynomial $P=P_m+\cdots+P_1$ as well as the convergence of level sets of $P_m$. We do not assume that $P$ is a polynomial solution until Theorem \ref{struct}. To start off, we write the irreducible factorization of $P_m$ as 
\begin{align}\label{factor}
	P_m=c_0p_{1}^{k_1}\cdots p_a^{k_a}
\end{align} 
where
\[\begin{aligned}
	&a\ge 1 \text{ is an integer},\\
	&c_0\neq0\text{ is} \text{ some constant},\\
	&k_i\ge 1 \text{ is an integer for }1\le i\le a,\\
	&p_i\text{ is an irreducible polynomial for } 1\le i\le a,\text{ and}\\
	&\gcd(p_i,p_j)=1 \text { for }1\le i\neq j\le a.
\end{aligned}\]
We group the factors $p_i^{k_i}$ as follows:
\begin{equation}\label{cond}
	\begin{aligned}
		I_1=&\{1\le i\le a:p_i\text{ changes sign}, k_i\text{ is odd}\} ,\\
		I_2=&\{1\le i\le a:p_i\text{ changes sign}, k_i\text{ is even}\},\\
		I_3=&\{1\le i\le a:p_i\text{ does not change sign}\}.
	\end{aligned}
\end{equation}
We can adjust the constant $c_0$ so that $p_i\ge0$ for all $i\in I_3$ which we assume in the sequel.
\begin{remark}
	Note that $I_j$ can be empty for $j=1,2,3$ and it will be shown in Theorem \ref{struct} that if $P$ is a polynomial solution, then $I_1$ consists of a single element and $I_2=\emptyset$.
\end{remark}

\begin{lemma}\label{reg}
	We have
	\begin{enumerate}
		\item[(i)]  For $i\in I_1\cup I_2 $, the set $\{p_i=0\}$ is $(n-1)$-dimensional and for $i\in I_3 $ the set $\{p_i=0\}$ is of dimension $\le n-2$.
		
		\item[(ii)]  For $ i,j\in I_1\cup I_2$, $i\neq j$, the set $\{p_i=0\}$ intersects $\{p_j=0\}$ along a set of dimension $\le n-2$.
	\end{enumerate}
\end{lemma}
\begin{proof}
	
	(i) This is a direct consequence of Lemma \ref{null} and (\ref{cond}).
	
	(ii) By Lemma \ref{null}, $\dim \{p_{i}=0\}= n-1$. If $p_{j}$ vanishes on an $(n-1)$-dimensional subset of $\{p_{i}=0\}$, then by Lemma \ref{null}, $p_{i}\mid p_{j}$. This is not possible since they are coprime by assumption (\ref{factor}). Then $\dim \{p_{i}=p_{j}=0\}< n-1$. Since $p_i$ and $p_j$ are polynomials, the set $\{p_{i}=p_{j}=0\}$ has dimension $\le n-2$.
\end{proof}

Now let $P=P_m+\cdots+P_1$ be a polynomial on $\R^n$ and $\graph P$ denote the graph of $P$ with the orientation given by the unit normal \[\frac{(\nabla P(x),-1)}{\sqrt{1+|\nabla P(x)|^2}}.\] We consider the following limit which we call \textit{translated blow-down}. Let $t\in\R$ be a number and $\lambda$ be a positive parameter that tends to 0. We consider the current $\eta_{\lambda\#}\graph P-t\lambda^{-m+1}e_{n+1}$ obtained by scaling down $\graph P$ by $\lambda$ and translating in the $-e_{n+1}$ direction by $t\lambda^{-m+1}$. We are interested in the limit as $\lambda\to0$. We have $\eta_{\lambda\#}\graph P-t\lambda^{-m+1}e_{n+1}=\graph G_t^{(\lambda)}$  where	
\begin{equation}\label{pt}
	\begin{aligned}
		G_t^{(\lambda)}(x)=&\lambda P(x/\lambda)-t\lambda^{-m+1}\\=&\lambda^{-m+1}(P_m(x)-t)+\lambda^{-m+2}P_{m-1}(x)+\cdots+P_1(x)
	\end{aligned}
\end{equation}
with the unit normal 
\begin{equation}\label{normal}
	\nu_{\graph G_t^{(\lambda)}}=\frac{(\nabla G_t^{(\lambda)}(x),-1)}{\sqrt{1+|\nabla G_t^{(\lambda)}(x)|^2}}.
\end{equation}
We also denote $G_0^{(\lambda)}$ by $P^{(\lambda)}$.  We have the following lemma on the support of a limit of $\graph G_t^{(\lambda_{j'})}$. 
\begin{lemma}\label{spt}
	Let $T\in I_{n}^\loc(\R^{n+1})$  be a limit of $[\graph G_t^{(\lambda_j)}]$ for some $\lambda_j\to0$. Then we have $\spt T\subset \{P_m=t\}\times \R$.
\end{lemma}
\begin{proof}
	A simple calculation shows
	\[
	\lim_{\lambda\to 0}|G_t^{(\lambda)}(x)|=\left\{
	\begin{aligned}
		&|P_1(x)| &\quad &\text{if }P_m(x)-t=P_{m-1}(x)=\cdots=P_2(x)=0\\
		&\infty &\quad &\text{otherwise}.
	\end{aligned}\right.\]
	Thus for any $x$ with $P_m(x)\neq t$, there is a neighborhood $U$ of $x$ in $\R^n$ so that $|G_t^{(\lambda)}|\to\infty$ uniformly in $U$ as $\lambda\to 0$. Thus  $(U\times \R)\cap \spt T=\emptyset$ and we have $\spt T\subset \{P_m=t\}\times\R$.
\end{proof}
Now we consider the convergence of $\graph G_t^{(\lambda)}$ as $\lambda\to0$.
\begin{lemma}\label{convgr}
Let $P=P_m+\cdots+P_1$ be a polynomial. Let $t\neq 0$ and $G_t^{(\lambda)}$ be as in (\ref{pt}). Then \[\lim_{\lambda\to0}\eta_{\lambda\#}\graph P-t\lambda^{-m+1}e_{n+1}=\lim_{\lambda\to0}\left[\graph G_t^{(\lambda)}\right] = \partial[P_m<t]\times \R\] where the convergence is both in $I_{n}^\loc(\R^{n+1})$ as currents and locally in $C^1$ as normal graphs over $\{P_m=t\}\times\R$.

Moreover, the family of hypersurfaces $(\graph G_t^{(\lambda)})_{\lambda\in [0,\epsilon)}$ is smooth in $\lambda$ and the velocity of the family at $\lambda=0$ is $-\frac{P_{m-1}}{|\nabla P_m|}$ viewed as a function on $\{P_m=t\}\times\R$.
\end{lemma}
\begin{proof}

Since $t\neq0$, by Lemma \ref{expansion} (i), we define $\nu(z)=\frac{\nabla P_m(z)}{|\nabla P_m(z)|}$ for all $z\in \{P_m=t\}$. By the implicit function theorem, we can choose a neighborhood $U$ of $z_0$ so that the following are satisfied:

(i) We can write $U$ as $U=\psi(W\times (-\sigma,\sigma))$ where $\sigma>0$, $W=U\cap\{P_m=t\}$ and $\psi:W\times (-\sigma,\sigma)\to U$ is a diffeomorphism given by $\psi(z,r)=z+r\nu(z)$.

(ii) There exists $\delta>0$ so that for all $(z,r)\in W\times (-\sigma,\sigma)$ we have \[\nabla P_m(z+r\nu(z))\cdot \nu(z)\ge 2\delta >0.\]
Using (\ref{pt}), we compute for $\lambda$ sufficiently small
\[
\begin{aligned}
	\frac{d}{dr}G_t^{(\lambda)}(z+r\nu(z))=&\nabla G_t^{(\lambda)}(z+r\nu(z))\cdot \nu(z)\\
	\ge&\lambda^{-m+1}\nabla P_m(z+r\nu(z))\cdot \nu(z)-O(\lambda^{-m+2})\\
	\ge&\delta \lambda^{-m+1}>0.
\end{aligned}
\]
Thus $r\mapsto G_t^{(\lambda)}(z+r\nu(z))$ is increasing for all $-\sigma<r<\sigma$.

For $z\in W$, we have \[|G_t^{(\lambda)}(z)|=|\lambda^{-m+2}P_{m-1}(z)+\cdots +\lambda^{-1}P_2(z)+P_1(z)|\le C\lambda^{-m+2}.\]
Therefore, the interval $[-\sigma \delta \lambda^{-m+1}+C\lambda^{-m+2},\sigma \delta \lambda^{-m+1}-C\lambda^{-m+2}]$ is in the range of $G_t^{(\lambda)}$ on $U$ by the fundamental theorem of calculus. Note that these intervals are increasing and cover $\R$ as $\lambda\to0$. Let $R>0$ and $(z,y)\in W\times (-R,R)$ be fixed.  There exists $\lambda_0(R)>0$ so that for all $0<\lambda<\lambda_0$, $(-R,R)$ is in the range of $r\mapsto G_t^{(\lambda)}(z+r\nu(z))$ on $(-\sigma,\sigma)$. Since $r\mapsto G_t^{(\lambda)}(z+r\nu(z))$ is increasing, there is a unique $-\sigma<r_0<\sigma$ so that $G_t^{(\lambda)}(z+r_0\nu(z))=y$.

We write $r_0=\vp^{(\lambda)}(z,y)$ and $x=x(z,y)=z+\vp^{(\lambda)}(z,y)\nu(z)$. Then
\begin{align}\label{graph}
	(x,G_t^{(\lambda)}(x))=(z,y)+\vp^{(\lambda)}(z,y)(\nu(z),0).
\end{align}
This means that we can write $\left[\graph G_t^{(\lambda)}\right]\llcorner (U\times(-R,R))$ as the normal graph of $\vp^{(\lambda)}$ over $W \times (-R,R)$. Using (\ref{graph}), the chain rule and (\ref{pt}), we compute the gradient of $\vp^{(\lambda)}$ as follows 
\begin{align}
	\frac{\partial \vp^{(\lambda)}}{\partial y}(z,y)=&\frac{1}{\nabla G_t^{(\lambda)}(x)\cdot\nu(z)}=O(\lambda ^{m-1})\\
	\frac{\partial \vp^{(\lambda)}}{\partial z_j}(z,y)=&\sum_{i=1}^{n}\frac{-1}{\nabla G_t^{(\lambda)}(x)\cdot\nu(z)}\frac{\partial G_t^{(\lambda)}}{\partial x_i}(x)(e_{ij}+\vp^{(\lambda)}(z,y)\frac{\partial \nu^i}{\partial z_j}(z))\\
	=&\sum_{i=1}^{n}\frac{-1}{\nabla P_m(x)\cdot\nu(z)}\frac{\partial P_m}{\partial x_i}(x)(e_{ij}+\vp^{(\lambda)}(z,y)\frac{\partial \nu^i}{\partial z_j}(z))+O(\lambda)\label{dz}
\end{align}
where  $(e_{ij})_{i=1,\dots,n,j=1,\dots,n-1}$ is the projection matrix from $\R^n$ to $T_z\{P_m=t\}=\nu(z)^{\perp}$.

From (\ref{pt}), we have 
\begin{equation}\label{defining.vp}
	P_m(x)-t=\lambda^{m-1}y-\lambda P_{m-1}(x)-\cdots- \lambda^{m-1}P_1(x).
\end{equation} Hence  with $x=z+\vp^{(\lambda)}(z,y)\nu(z)$, we have $P_m(x)\to t$ as $\lambda\to0$. Since $\nabla P_m\neq0$ on $\{P_m=t\}$, we have $d(x,\{P_m=t\})\to 0$. Hence\[\lim_{\lambda\to 0}\vp^{(\lambda)}(z,y)= 0.\]
Since $(e_{ij})$ is the projection onto $T_z\{P_m=t\}$ we have $\sum_{i=1}^{n}e_{ij}\partial_{x_i}P_m(z)=0$. Thus the first term in (\ref{dz}) converges to 0 in view of $d(x,\{P_m=t\})\to 0$. Hence we see that $\graph G_t^{(\lambda)}$ converge to $(U\cap \{P_m=t\})\times\R$ as normal graphs in $C^1$. We observe that the unit normal of $\graph G_t^{(\lambda)}$ as in (\ref{normal}) converges to $\nu$ as $\lambda\to0$. Thus $\left[\graph G_t^{(\lambda)}\right]\to \partial[P_m<t]\times \R$ as currents.

To show that the velocity of the family $(\graph G_t^{(\lambda)})_{\lambda\in [0,\epsilon)}$ at $\lambda=0$ is $-\frac{P_{m-1}}{|\nabla P_m|}$, it suffices to compute $\frac{\partial \vp^{(\lambda)}}{\partial \lambda}|_{\lambda=0}$. By Taylor expansion around $z\in \{P_m=t\}$, we have
\[\begin{aligned}
	P_m(x)-t=&P_m(z+\vp^{(\lambda)}\nu(z))-P_m(z)\\
	=&\nabla P_m(z)\cdot\nu(z)\vp^{(\lambda)}+o(\vp^{(\lambda)})\\
	=&|\nabla P_m(z)|\vp^{(\lambda)}+o(\vp^{(\lambda)}).
\end{aligned}\]
Since right hand side of (\ref{defining.vp}) is of order $O(\lambda)$ and $|\nabla P_m(z)|\neq0$, we have $\vp^{(\lambda)}=O(\lambda)$.
We divide both sides of (\ref{defining.vp}) by $|\nabla P_m(z)|\lambda$ and obtain
\[\frac{\vp^{(\lambda)}}{\lambda}+o(1)=-\frac{P_{m-1}(x)}{|\nabla P_m(z)|}+O(\lambda).\]
Taking $\lambda\to0$ we get $\frac{\partial \vp^{(\lambda)}}{\partial \lambda}|_{\lambda=0}(z)=-\frac{P_{m-1}}{|\nabla P_m|}$. This finishes the proof.
\end{proof}

\begin{definition}\label{regular.point}
 	We say that $z_0\in \{P_m=0\}$ is a \textit{regular point} of $\{P_m=0\}$ if for some $ i\in I_1\cup I_2$
	\begin{equation}
		\label{regular}\begin{aligned}
			&p_i(z_0)=0,\ \nabla p_i(z_0)\neq 0,\text{ and}\\
			&p_j(z_0)\neq 0\text{ for all }1\le j\le a,\ j\neq i.
		\end{aligned} 
	\end{equation}
	A point $z_0\in \{P_m=0\}$ is a \textit{singular point} if it is not a regular point.
\end{definition}
\begin{remark}
	 A regular point in the sense of Definition \ref{regular.point} is a regular point of the current $\partial [P_m<0]$ in the sense of \cite[5.3.17]{Fed}. The converse is not true since there could be $j\in I_3$ such that $p_j(z_0)=0$ for $z_0\in \reg \partial [P_m<0]$. It will be shown in Theorem \ref{struct} (ii) (c) that if $P$ is a polynomial solution, then the regular points of $\{P_m=0\}$ in the sense of Definition \ref{regular.point} coincides with the regular points of $\partial[P_m<0]$ as a current.
\end{remark}

Now we consider the convergence of $\partial[P_m<t]$ as $t\to0$ in a neighborhood of $z_0$ where $z_0$ is a regular point of $\{P_m=0\}$.
In view of the factorization (\ref{factor}), we denote the polynomial $c_0p_{1}^{k_1}\cdots  \widehat{p_i^{k_i}}\cdots  p_{a}^{k_{a}}$ by $P_m/p_i^{k_i}$ for $i\in I_1\cup I_2$.

\begin{lemma}\label{convls}Let $z_0$ be as in (\ref{regular}). Then there exists a neighborhood $U$ of $z_0$ in $\R^n$ so that the following items hold.
	\begin{enumerate}
		\item [(i)] If $i\in I_1$, then
		\[\lim_{t\to0}\left(\partial[P_m<t]\right)\llcorner{U}= \sgn \left(\left(P_m/p_i^{k_i}\right)(z_0)\right)\left(\partial[p_i<0]\right)\llcorner{U}\]
		where the convergence is both in $I_{n-1}^\loc(U)$ as currents and in $C^{\infty}$ as smooth manifolds.
		
		In particular
		\[\lim_{t\to0}\|\partial[P_m<t]\|\llcorner U= \HH^{n-1}\llcorner(U\cap \{p_i=0\}).\]
		\item[(ii)]  If $i\in I_2$  then as currents in $I_{n-1}^\loc(U)$ we have
		\[\lim_{t\to0}\partial[P_m<t]\llcorner U=0.\]
		(a) If $\sgn ((P_m/p_i^{k_i})(z_0))>0$, then
		\[
		\begin{aligned}
			\lim_{t\to0^+}\|\partial[P_m<t]\|\llcorner U=& 2 \HH^{n-1}\llcorner (U\cap \{p_i=0\}),\\
			\lim_{t\to0^-}\|\partial[P_m<t]\|\llcorner U=&0.
		\end{aligned}
		\]
		(b) If $\sgn ((P_m/p_i^{k_i})(z_0))<0$, then
		\[
		\begin{aligned}
			\lim_{t\to0^+}\|\partial[P_m<t]\|\llcorner U=& 0,\\
			\lim_{t\to0^-}\|\partial[P_m<t]\|\llcorner U=& 2 \HH^{n-1}\llcorner (U\cap \{p_i=0\}).
		\end{aligned}
		\]
	\end{enumerate}
\end{lemma}

\begin{proof}
	(i)  By (\ref{cond}), $k_i$ is odd for $i\in I_1$. We have $\{P_m=t\}=\{w=t^{1/k_i}\}$ where 
	\begin{equation}\label{w}
		w=\left(P_m/p_i^{k_i}\right)^{1/k_i}p_i.
	\end{equation}
	The function $w$ is smooth around $z_0$, since $P_m/p_i^{k_i}$ is nonzero around $z_0$ and $k_i$ is odd. Since $\nabla p_i(z_0)\neq 0$ and $p_j(z_0)\neq0$ for $j\neq i$, we have $\nabla w(z_0)\neq0$. By the implicit function theorem, there exist $\delta>0$, a neighborhood $U$ of $z_0$, and $f:(U\cap \{w=0\})\times (-\delta,\delta)\to U$ a diffeomorphism so that $w(f(z,s))=s$ and $f(z,0)=z$ for all $(z,s)\in (U\cap \{w=0\})\times(-\delta,\delta)$. Therefore, $U\cap \{w=t^{1/k_i}\}=f(U\cap \{w=0\},t^{1/k_i})$. Since $f(\cdot,t)\to \id_{U\cap \{w=0\}}$ in $C^{\infty}$ as $t\to0$, we have \[U\cap\{P_m=t\}=U\cap\{w=t^{1/k_i}\}\to U\cap\{w=0\}=U\cap \{p_i=0\}\quad\text{in } C^{\infty}\text{ as }t\to0.\]
	
	Note that on $\{P_m=t\}$ for $t\neq 0$, we have \[\frac{\nabla P_m}{|\nabla P_m|}=\frac{\nabla (w^{k_i})}{|\nabla (w^{k_i})|}=\frac{k_iw^{k_i-1}\nabla w}{|k_iw^{k_i-1}\nabla w|}=\frac{\nabla w}{|\nabla w|}\] where the last equality is regardless of the sign of $w$ since $k_i$ is odd (hence $k_i-1$ is even). Thus we have \[\frac{\nabla P_m}{|\nabla P_m|}=\frac{\nabla w}{|\nabla w|}\to \sgn \left(\left(P_m/p_i^{k_i}\right)(z_0)\right)\frac{\nabla p_i}{|\nabla p_i|}\text{ as }t\to0.\] This establishes the convergence with the desired orientation.
	
	(ii) We only prove (a). Item (b) follows from applying (a) to $-P_m$.
	
	Since $\sgn ((P_m/p_i^{k_i})(z_0))>0$ and $k_i$ is even, we have $U\cap\{P_m=t\}=\emptyset$ for $t<0$ for a small neighborhood $U$ of $z_0$. Thus the convergence for $t\to0^-$ is obtained. 
	
	For $t\to0^+$, we consider $w=\left(P_m/p_i^{k_i}\right)^{1/k_i}p_i$. This function is well defined and smooth around $z_0$ since $\sgn ((P_m/p_i^{k_i})(z_0))>0$. As in (i), we choose a neighborhood $U$ of $z_0$ so that $U\cap \{w=s\}$ is the image of $U\cap \{w=0\}$ under a smooth diffeomorphism $f(\cdot,s)$ for $s\in (-\delta,\delta)$. Now for $t> 0$ small, we have  \[U\cap\{P_m=t\}=(U\cap\{w=t^{1/k_i}\})\cup(U\cap\{w=-t^{1/k_i}\}).\]
	
	On $\{w=s\}$ we have \[\frac{\nabla P_m}{|\nabla P_m|}=\frac{\nabla (w^{k_i})}{|\nabla (w^{k_i})|}=\frac{k_iw^{k_i-1}\nabla w}{|k_iw^{k_i-1}\nabla w|}=\sgn(w)\frac{\nabla w}{|\nabla w|}\] where the last equality is due to $k_i$ is even (hence $k_i-1$ is odd). Since $\sgn ((P_m/p_i^{k_i})(z_0))>0$, we have $\frac{\nabla w}{|\nabla w|}\to \frac{\nabla p_i}{|\nabla p_i|}$ as $s\to0$. Thus by the same argument as in (i), we have \[U\cap\{w=\pm t^{1/k_i}\}\to \pm U\cap \{p_i=0\}\quad\text{ as }t\to0^+\] both in $I_{n-1}^\loc(U)$ and in $C^\infty$ as manifolds. Thus we have $\partial[P_m<t]\llcorner U\to0$ since opposite orientations cancel with each other and that $\|\partial[P_m<t]\|\llcorner U\to 2 \HH^{n-1}|_{U\cap \{p_i=0\}}$ as $t\to0^+$.
\end{proof}

\begin{lemma}\label{partial[P_m<0]and|partial[P_m<0]|}
If $P_m$ is as in \eqref{factor} and \eqref{cond}, then
    \[\partial[P_m<0]=\sum_{i\in I_1}\sgn \left(P_m/p_i^{k_i}\right)\partial [p_i<0]\]
    and
    \[\|\partial[P_m<0]\|=\sum_{i\in I_1}\|\partial [p_i<0]\|.\]
\end{lemma}
\begin{proof}
    By Lemma \ref{reg}, the set of singular points of $\{P_m=0\}$ has dimension $\le n-2$. We only need to check the identities in a small neighborhood $U$ of any regular point $z_0\in \{p_i=0\}$ of $\{P_m=0\}$.
	
	If $i\in I_1$, we define $w$ in $U$ as in (\ref{w}). We have $[P_m<0]\llcorner U=\sgn ((P_m/p_i^{k_i})(z_0)) [w<0]\llcorner U$. Then we have $\partial [P_m<0]\llcorner U=\sgn ((P_m/p_i^{k_i})(z_0))\partial [p_i<0]\llcorner U$ and $\|\partial [P_m<0]\|\llcorner U=\|\partial [p_i<0]\|\llcorner U$.
	
	If $i\in I_2$, and $\sgn ((P_m/p_i^{k_i})(z_0))>0$, we have $P_m\ge0$ in $U$ and thus $[P_m<0]\llcorner U=0$. Thus $\partial [P_m<0]=0$ and $\|\partial [P_m<0]\|\llcorner U=0$.
	
	If $i\in I_2$, and $\sgn ((P_m/p_i^{k_i})(z_0))<0$, we have $P_m\le0$ in $U$. Then we see that $[P_m<0]\llcorner U=[P_m\neq0]\llcorner U=[U]$. Thus $\partial [P_m<0]=0$ and $\|\partial [P_m<0]\|\llcorner U=0$.
\end{proof}
\begin{corollary}\label{convlscor}We have the following convergences
	\begin{enumerate}
		\item[(i)]
		\[
		\lim_{t\to 0}\partial[P_m<t] =\sum_{i\in I_1}\sgn \left(P_m/p_i^{k_i}\right)\partial [p_i<0].
		\]
		\item[(ii)]
		\[\lim_{t\to 0+}\|\partial[P_m<t]\|= \sum_{i\in I_1}\HH^{n-1}\llcorner{\{p_i=0\}}+\sum_{i\in I_2}2\HH^{n-1}\llcorner{\{p_i=0,P_m/p_i^{k_i}>0\}}.\]
		\item[(iii)]
		\[\lim_{t\to 0-}\|\partial[P_m<t]\|= \sum_{i\in I_1}\HH^{n-1}\llcorner{\{p_i=0\}}+\sum_{i\in I_2}2\HH^{n-1}\llcorner{\{p_i=0,P_m/p_i^{k_i}<0\}}.\]
	\end{enumerate}
\end{corollary}
\begin{proof}
	For any $z\in\{P_m\neq0\}$, the set $\{P_m=t\}$ has positive distance to $z$ for $t$ sufficiently small. Thus for any limit $T\in I_{n-1}^\loc(\R^{n})$ of $\partial [P_m<t_i]$ as $t_i\to0$, we have $\spt T\subset \{P_m=0\}$. By Lemma \ref{reg}, the set of singular points of $\{P_m=0\}$ has dimension $\le n-2$. Thus it suffices to consider the limits only near regular points of $\{P_m=0\}$. Then by Lemma \ref{convls}, all three limits in the statement exist and are equal to the desired ones by summing over $i\in I_1\cup I_2$ in Lemma \ref{convls}. This finishes the proof.
\end{proof}

\subsection{The structure theorem}
Now we can prove the main structure theorem.
\begin{theorem}\label{struct}
	Let $P=P_m+\cdots+P_1$ be a polynomial solution. Then we have:
	\begin{enumerate}
		\item[(i)] For every $t\in\R$, the current $\partial [P_m<t]$ is area-minimizing. If $t\neq0$, $\partial[P_m<t]$ is a smooth minimal hypersurface.
		\item[(ii)] There exist 
		\begin{enumerate}
			
			\item [(a)]$k\ge 1$ odd, 
			
			\item[(b)]$p $ an irreducible homogeneous polynomial that changes sign and has $\deg p\ge2$,
			
			\item[(c)] $Q_m$ an even degree homogeneous polynomial which is nonnegative on $\R^n$ and coprime with $p $, with $\{Q_m=0\}\subset\{p=0\}\cap\{\nabla p=0\}$ of dimension $\le n-8$,
			
		\end{enumerate}
		so that
		\begin{align}\label{pkq}
			P_m=p ^kQ_m.
		\end{align}
		
		\item[(iii)] The unique tangent cone of $\graph P$ at infinity is $\partial [p <0]\times \R$.
		
		\item[(iv)]The function $u=\frac{1}{|\nabla (p Q_m^{1/k})|}$ is the velocity of the family of minimal hypersurfaces $(\{P_m=t^k\})_{t\in\R}$ and is finite on the regular part of $\{P_m=t^k\}$ for every $t\in \R$.

		\item[(v)]  If $k\ge 3$ in (\ref{pkq}), then there exists $2\le s\le m-1$ so that
		\[P=p ^kQ_m+p ^{k_{m-1}} Q_{m-1}+\cdots +p ^{k_{s+1}}Q_{s+1}+p Q_s+P_{s-1}+\cdots+P_1\]
		where  $p $ and $Q_i$ are coprime for $s\le i\le m-1$ and $k_i\ge 2$ for $s+1\le i\le m-1$.
	\end{enumerate}
\end{theorem}
\begin{remark}
	We remark that in the introduction, Theorem \ref{struct0} states $\deg p\ge 3$ instead of $\deg p\ge2$ here. This will be proved later in Theorem \ref{noisopara}.
\end{remark}
\begin{remark}
	(ii) (c) tells us that the regular part of $\{P_m=0\}$ in the sense of Definition \ref{regular.point} is the same as the regular part of $\partial [P_m<0]$ as a current.
\end{remark}
\begin{proof}
	(i) Since $P$ is a polynomial solution, we have $\graph P$ as well as its translations and rescalings are area-minimizing boundaries in $\R^{n+1}$ by Lemma \ref{minimizing}. By Lemma \ref{convgr} and Lemma \ref{FedFlem}, we have for any $t\neq 0$, $\partial [ P_m<t]\times\R$ is area-minimizing. Hence $\partial [ P_m<t]$ is area-minimizing for $t\neq0$. By Corollary \ref{convlscor} and Lemma \ref{partial[P_m<0]and|partial[P_m<0]|}, we have $\partial [ P_m<t]\to \partial [ P_m<0]$ as $t\to0$. Thus by Lemma \ref{FedFlem} again, $\partial [ P_m<0]$ is area-minimizing. The smoothness of nonzero level sets of $P_m$ is obtained in Lemma \ref{expansion} (ii).
	
	(ii)  We first claim that $I_1\cup I_2\neq\emptyset$. Indeed, if this is not true, then by assumption (\ref{cond}), we can assume $P_m\ge0$. By Corollary \ref{convlscor} (ii), (iii) we have $\lim_{t\to0}\|\partial [ P_m<t]\|=0$. However since $P_m$ is homogeneous of degree $m$ we have $\{P_m=t\}=\eta_{t^{1/m}}(\{P_m=1\})$ for $t>0$. Then \[\|\partial [ P_m<t]\|(B_1)=t^{(n-1)/m}\|\partial [ P_m<1]\|(B_{t^{-1/m}})\ge c(n)>0\] where in the last inequality we use the monotonicity formula for the minimal hypersurface $\partial [P_m<1]$ \cite[Theorem 17.6]{Sim}. This leads to a contradiction. 
	
	By Lemma \ref{FedFlem}, Lemma \ref{partial[P_m<0]and|partial[P_m<0]|}  and Corollary \ref{convlscor} (i), we have as $t\to0$ \[\|\partial [ P_m<t]\|\to \|\partial [ P_m<0]\|=\sum_{i\in I_1}\HH^{n-1}\llcorner{\{p_i=0\}}.\] In view of Corollary \ref{convlscor} (ii), (iii), we see that $I_2=\emptyset$. In other words, there is no factor $p_i^{k_i}$ where $p_i$ changes sign and $k_i$ is even. Thus $I_1\neq\emptyset$.
	
	We claim that $I_1$ consists of a single element. Indeed, we show that any polynomial, not necessarily the leading term of a polynomial solution, of the form $P_m=c_0\prod_{i\in I_1\cup I_3}p_i^{k_i}$ such that $I_1,I_3$ are as in \eqref{cond} and $\partial [P_m<0]$ is area-minimizing must satisfy  $|I_1|=1$ where $|I_1|$ is the number of elements in $I_1$. We note that $\partial[P_m>0]=-\partial[P_m<0]$ is also area-minimizing. We use induction to show that $|I_1|\neq k$ for any $\N\ni k\ge 2$. We first show that $|I_1|\neq 2$. Suppose for contradiction that $I_1=\{1,2\}$. We write $\tilde p=c_0\prod_{i\in I_1\cup I_3,i\neq 1}p_i^{k_i}$. We have \[\{P_m<0\}=\{p_1<0,\tilde p>0\}\cup \{p_1>0,\tilde p<0\}.\]  Since $\partial [P_m<0]$ is area-minimizing, Lemma \ref{indecomposable} implies that $\{P_m<0\}$ is connected. Hence either $\{p_1<0,\tilde p>0\}=\emptyset$ or $\{p_1>0,\tilde p<0\}=\emptyset$. Let's assume $\{p_1<0,\tilde p>0\}=\emptyset$. The other case can be handled similarly. Then we have $\{p_1<0\}\subset \{\tilde p\le0\}$ and $\{\tilde p>0\}\subset \{p_1\ge0\}$. By continuity, we have $\{p_1\le 0\}\subset \{\tilde p\le0\}$ and $\{\tilde p\ge 0\}\subset \{p_1\ge0\}$. Thus taking complement, we have $\{\tilde p>0\}\subset \{p_1>0\}$ and $\{p_1<0\}\subset \{\tilde p<0\}$. Thus \begin{equation}\label{P_m>0}
	    \{P_m>0\}=\{p_1>0,\tilde p>0\}\cup \{p_1<0,\tilde p<0\}=\{\tilde p>0\}\cup\{ p_1<0\}.
	\end{equation} Observe that the two open sets on the right hand side are disjoint. Since $I_1=\{1,2\} $, $\tilde p=c_0p_2^{k_2}\prod_{i\in I_3}p_i^{k_i}$. By Lemma \ref{null}, $p_i^{k_i}$ does not change sign and has zero locus of codimension at least 2 for $i\in I_3$ and $p_1$, $p_2^{k_2}$ changes sign and has zero locus of dimension $n-1$. Thus we have $\{\tilde p>0\}\neq\emptyset\neq \{p_1<0\}$.  This contradicts Lemma \ref{indecomposable} for the area-minimizing boundary $-\partial[P_m<0]=\partial[P_m>0]$ and finishes the proof of $|I_1|\neq 2$. Suppose we know $|I_1|\neq k-1$ for some $k\ge 3$. To show that $|I_1|\neq k$, we repeat the above argument up to \eqref{P_m>0}. By Lemma \ref{partial[P_m<0]and|partial[P_m<0]|}, $\partial[P_m>0]=\partial [\tilde p>0]+\sgn (P_m/p_i^{k_i})\partial [p_1>0]$ and $\|\partial[P_m>0]\|=\|\partial [\tilde p>0]\|+\|\partial [p_1>0]\|$. Since $\partial [P_m>0]$ is area-minimizing, $\partial [\tilde p>0]$ is also area-minimizing by \cite[Lemma 33.4]{Sim}. For $\tilde p$, we have $|I_1|=k-1$ which is a contradiction to the induction hypothesis. This finishes the induction.

	Thus $I_1$ consists of a single element which we call $1$. We let $p=p_1$,  $k=k_1$ and $Q_m=c_0\prod_{i\in I_3}p_{i}^{k_{i}}$. We can assume $c_0>0$ in (\ref{factor}) by adjusting $p_1$ to be $-p_1$ if necessary. Then $Q_m\ge 0$ and hence $\deg Q_m$ is even. We get $P_m=p^kQ_m$.
	
	Now we show that $\{Q_m=0\}\subset \{p=0\}$. Suppose for contradiction that there is $x_0\in \{Q_m=0\}$ so that $p(x_0)\neq0$. Without loss of generality, we assume $p(x_0)>0$. By continuity of $p$, we choose $r>0$ small so that for any $x\in B_{2r}(x_0)$ we have $p(x)\ge \frac{1}{2}p(x_0)$. Since $\{Q_m=0\}$ has dimension $\le (n-2)$, by continuity of $P_m$, for any $t>0$ small, there exists $y\in B_{r/2}(x_0)$ so that $P_m(y)=t$. By the monotonicity formula for the minimal hypersurface $\partial [P_m<t]$ \cite[Theorem 17.6]{Sim} we have \[\|\partial [P_m<t]\|(B_r(x_0))\ge \|\partial [P_m<t]\|(B_{r/2}(y))\ge c(n)r^n.\]
	By Lemma \ref{FedFlem} and the first equality of Corollary \ref{convlscor} (i), we have
	\[\|\partial [P_m<0]\|(B_{2r}(x_0))\ge \limsup_{t\to0}\|\partial [P_m<t]\|(B_r(x_0))\ge c(n)r^n>0.\]
	However, by the second equality of Corollary \ref{convlscor} (i), we have $\spt \|\partial [P_m<0]\|=\spt \|\partial[p<0]\|\subset \{p=0\}$. This is a contradiction.

	All desired results have been established except for $\deg p \ge 2$, which will be proved after (iii), and $\{Q_m = 0\} \subset \{\nabla p = 0\}$ and $\dim \{Q_m=0\}\le n-8$, which will be proved after (iv).
	
	(iii) Let $T$ be a tangent cone at infinity of $\graph P$. By Lemma \ref{spt}, $\spt T\subset \{P_m=0\}\times \R$. By \cite[\S 2]{simon1987asymptotic} and Lemma \ref{indecomposable} we know that $T=\partial [U]\times \R$ for some connected open set $U\subset\R^n$. Hence by (\ref{factor}) and the results in (ii) we have \[\spt \partial [U]\subset \{P_m=0\}=\{p =0\}.\] If we take $\phi\in C_c^\infty(\{p <0\})$ then $\partial [U](\phi)=0$. Thus $\partial[U]=0$ in $I_{n-1}^\loc(\{p <0\})$. Similarly $\partial [U]=0$ in $I_{n-1}^\loc(\{p >0\})$. Thus by the constancy theorem \cite[Theorem 26.27]{Sim}, $[U]$ is constant on both $\{p <0\}$ and $\{p >0\}$. Hence either $[U]=[p <0]$ or $[U]=[p >0]$. Recall that we wrote $\eta_{\lambda\#}\graph P= \graph P^{(\lambda)}$ in (\ref{pt}) with the unit normal as in (\ref{normal}). For $x\notin \{p =0\}$ which is close to $z\in \{p=0\}\cap\{\nabla p\neq0\}$, $\nu_{\graph P^{(\lambda)}}(x)$ is close to  $\frac{\nabla p(z) }{|\nabla p(z) |}$ for $\lambda$ small. Since $\frac{\nabla p}{|\nabla p|}$ is the outward unit normal of $\{p<0\}$, we have $U\cap \{p <0\}\neq \emptyset$. Therefore $ [U]=[p <0]$.
	
	To show $\deg p\ge 2$ in (ii) (b), we assume for contradiction that $\deg p=1$. Then the tangent cone at infinity of $\graph P$ is $\partial[p <0]\times \R$ which is a hyperplane. Then the monotonicity formula \cite[Theorem 17.6]{Sim} shows that $\graph P$ is a hyperplane. Hence $\deg P=1$ which contradicts Definition \ref{polynomial.solution}.
	
	(iv) 
	By Theorem \ref{struct} (i) for any $t\in\R$, the level set $\{p Q_m^{1/k}=t\}=\{P_m=t^k\}$ is minimal. We take a $\frac{1}{k}$ power of $P_m$ since $\nabla (p Q_m^{1/k})=Q_m^{1/k}\nabla p \neq 0$ on the regular part of $\{P_m =0\}$. Now the velocity of the family of minimal hypersurfaces $(\{p Q_m^{1/k}=t\})_{t\in\R}$ is $u=\frac{1}{|\nabla (p Q_m^{1/k})|}$ which is finite on the regular part of $\{p Q_m^{1/k}=t\}$ for any $t\in \R$.
	
	Now we show that $\{Q_m=0\}\subset \{\nabla p=0\}$ in (ii) (c). Suppose for contradiction that there is $x_0\in \{Q_m=0\}$ such that $\nabla p(x_0)\neq0$. For $t\in\R$, we write $\Sigma_t=\{P_m=t^k\}$. Since $\{Q_m=0\}\subset\{p=0\}$, we have $x_0$ is a regular point of $\partial[p<0]$ as a current and $\partial[P_m<t^k]\rightharpoonup\partial[p<0]$ as $t\to0$ by Corollary \ref{convlscor} (i). By Allard's regularity \cite[Theorem 24.2]{Sim} and elliptic regularity \cite{GT}, there exist $r>0$ small and $f_t:\Sigma_0\cap B_r(x_0)\to \R$ such that $\Sigma_t\cap B_r(x_0)=\graph_{\Sigma_0}f_t$ and $f_t\to0$ smoothly as $t\to0$. We restrict ourself to $t\to 0^+$ and hence $\Sigma_t$ is on one side of $\Sigma_0$. We have $f_t>0$. By what we have just proved in (iv), we have $\frac{d}{dt}|_{t=0}f_t=\lim_{t\to0^+} f_t/t=u=\frac{1}{Q_m^{1/k}|\nabla p|}$. 
	
	For any $t>0$ small, since $\Sigma_t$ is minimal, $f=f_t$ satisfies an equation of the form \[\sum_{i,j=1}^{n-1}(a_{ij}(\nabla f,f,x)f_{x_i})_{x_j}+\sum_{i=1}^{n-1}b_i(\nabla f,f,x)f_{x_i}+c(f,x)f=0\] in local coordinates of $\Sigma_0$. There exist $\Lambda>1,\epsilon>0$ such that if $|f|_{C^1}<\epsilon$ then $|a_{ij}|+|b_i|+|c|\le \Lambda$ and $\Lambda|\xi|^2\ge a_{ij}\xi_i\xi_j\ge \Lambda^{-1}|\xi|^2$. Since $f_t\to0$ smoothly on $\Sigma_0\cap B_{r}(x_0)$, the above conditions are satisfied for $t$ sufficiently small.
	By the Harnack inequality \cite[Theorem 8.20]{GT}, there exists $C>0$ independent of $t$ such that 
	\[\sup_{\Sigma_0\cap B_{r/2}(x_0)}f_t\le C\inf_{\Sigma_0\cap B_{r/2}(x_0)} f_t.\] Thus we have
	\[\sup_{\Sigma_0\cap B_{r/2}(x_0)}t^{-1}f_t\le C\inf_{\Sigma_0\cap B_{r/2}(x_0)}t^{-1}f_t.\]Taking $t\to 0$, we obtain \[\sup_{\Sigma_0\cap B_{r/2}(x_0)} u\le C\inf_{\Sigma_0\cap B_{r/2}(x_0)}u<\infty.\]
	In particular, we have $u<\infty$ in $\Sigma_0\cap B_{r/2}(x_0)$. This contradicts $u(x_0)=\frac{1}{Q_m^{1/k}(x_0)|\nabla p(x_0)|}=\infty$.
	
	With $\{Q_m=0\}\subset \{p=\nabla p=0\}$, we have $\dim\{Q_m=0\}\le n-8$ since $\partial [p<0]$ is an area-minimizing boundary and hence its singular set $\{p=\nabla p=0\}$ is of dimension $\le n-8$ by \cite{Fed70}.
	
	(v)
	Suppose $k\ge 3$ in (\ref{pkq}). Then $\nabla P_m\equiv0$ on $\{p =0\}$.  By (iii) we have $[\graph P^{(\lambda)}]\to \partial [p <0]\times \R$ as $\lambda\to 0$. By Allard's regularity \cite[Theorem 24.2]{Sim}, for any regular point $(z,y)\in \{p =0\}\times \R$, there is a neighborhood $U$ of $(z,y)$ so that $\graph P^{(\lambda)}\cap U$ can be written as the normal graph of $\vp^{(\lambda)}$ over $(\{p =0\}\times\R)\cap U$ and $\vp^{(\lambda)}\to0$ in $C^{1,\alpha}$ as $\lambda\to0$. In particular, with $x=z+\vp^{(\lambda)}(z,y)\nu(z)$, we have
	\[
	\frac{(\nabla P^{(\lambda)}(x),-1)}{\sqrt{1+|\nabla P^{(\lambda)}(x)|^2}}=\nu_{\graph P^{(\lambda)}}(x)\to \nu(z)=\frac{(\nabla p (z),0)}{|\nabla p (z)|}\text{ as }  \lambda\to 0.
	\]
	Since $\frac{1}{\sqrt{1+|\nabla P^{(\lambda)}|^2}}\to 0$, $\nabla P_m(z)=0$ and $|\nabla P_1|$ is homogeneous of degree 0, there exists $2\le s\le m-1$ so that 
	\begin{equation}\label{grad}
		\nabla P_m(z)=\nabla P_{m-1}(z)=\cdots=\nabla P_{s+1}(z)=0\text{ and }\nabla P_s(z)\neq 0.
	\end{equation} 
	
	Since for every regular point $z\in \{p =0\}$, there exists $s$ satisfying (\ref {grad}). We have
	\[\{p=0\}\cap \{\nabla p\neq0\}=\bigcup_{s=2}^{m-1}V_s\] where \[V_s:=\{z\in \{p=0\}\cap \{\nabla p\neq0\}:\text{(\ref{grad}) holds}\}.\] Since $\dim \{p=0\}\cap\{\nabla p\neq0\}=n-1$ and $\{V_s\}_{s=2}^{m-1}$ are disjoint, we may choose the maximal $s$ so that the semi-algebraic set $V:=V_s$ has dimension $(n-1)$ i.e. (\ref{grad}) holds on a semi-algebraic set of dimension $(n-1)$. By Lemma \ref{null}, we have $p \mid \nabla P_i$ for $s+1\le i\le m-1$. By Euler's formula $iP_i(x)=x\cdot\nabla P_i(x)$, we have $p \mid P_i$. Thus we can write $P_i=p \tilde Q_i$. Since \[0=\nabla P_i=\tilde Q_i\nabla p  +p \nabla \tilde Q_i=\tilde Q_i\nabla p \quad \text{in }V,\]we see that $\tilde Q_i$ vanishes on $V$. Thus by Lemma \ref{null}, $p \mid \tilde Q_i$ and hence $p ^2\mid P_i$.
	
	By  (\ref{grad}), we have for $(z,y)\in V\times \R$, as $\lambda\to0$ 
	\[\frac{(\nabla P^{(\lambda)}(x),-1)}{\sqrt{1+|\nabla P^{(\lambda)}(x)|^2}}\to \frac{(\nabla P_s(z),0)}{|\nabla P_s(z)|}.\]
	Hence for all $z\in V$, we have
	\[\frac{\nabla P_s(z)}{|\nabla P_s(z)|}=\frac{\nabla p (z)}{|\nabla p (z)|}.\]
	By Euler's identity, \[sP_s(z)=z\cdot \nabla P_s(z)=z\cdot \nabla p (z)\frac{|\nabla P_s(z)|}{|\nabla p (z)|}=\deg p \cdot p (z)\frac{|\nabla P_s(z)|}{|\nabla p (z)|}=0.\] Thus we have $P_s(z)=0$ on $V$. Thus by Lemma \ref{null}, $p \mid P_s$. Since $\nabla P_s\neq0$ on $V$, we have $p^2\nmid P_s$. This finishes the proof.
\end{proof}

\section{Degree estimates of polynomial solutions}
\label{4}
In this section, we use Theorem \ref{struct} to give estimates for the degree of polynomial solutions.

We first show that cubic polynomial solutions do not exist. This is important in later applications (cf. Theorem \ref{R8cubic}).
\begin{theorem}\label{nodeg23}
	Let $P$ be a polynomial solution. Then $\deg P\ge 4$.
\end{theorem}

\begin{proof}By Definition \ref{polynomial.solution} and Theorem \ref{nodeg2}, we only need to rule out cubic polynomial solutions. Suppose for the purpose of contradiction that $P=P_3+P_2+P_1$ is a cubic polynomial solution. By Theorem \ref{struct} (ii), we know that $P_3=p ^kQ_3$ where $Q_3\ge0$ has even degree, $\deg p\ge 2$ and $k\ge 1$ is odd. Thus we must have $\deg Q_3=0$, $k=1$ and $P_3$ is irreducible. We have $LP_3=0$ by Lemma \ref{expansion} (i). By  \cite[Theorem 1]{tkachev2010classification} (see also \cite[Theorem 6.6.1]{NTV}), we know that any irreducible cubic polynomial solving $LP_3=\lambda|x|^2P_3$ for some $\lambda\in\R$ is harmonic. Thus $\Delta P_3=0$.  Hence $P_3=0$ by Lemma \ref{infhar}, a contradiction.
\end{proof}
Next we use Theorem \ref{struct} (iv) to bound the degree of $p$.

\begin{theorem}\label{n-2}
	Let $P$ be a polynomial solution  on $\R^n$ and $P_m=p^kQ_m$ as in Theorem \ref{struct} (ii). Then 
	\[\mu_n^-< \deg p +k^{-1}\deg Q_m< \mu_n^+\]where
	$\mu_n^{\pm}= \frac{n-1\pm\sqrt{(n-3)^2-4(n-2)}}{2}$. In other words, $k\mu_n^-<\deg P< k\mu_n^+$.
\end{theorem}
\begin{proof}
	Let $\Sigma=\{p=0\}\cap\{\nabla p\neq0\}\cap \SSS^{n-1}$ be the regular part of the link of the cone $\{p=0\}$. By Theorem \ref{struct} (iv), $u=\frac{1}{Q_m^{1/k}|\nabla p|}$ is a Jacobi field on regular part of $\{p=0\}$ i.e. $\Delta_{\{p=0\}}u+|A_{\{p=0\}}|^2u=0$. We write $u=r^{-\gamma}\phi$ where \[r=|x|,\quad \gamma=\deg p+\frac{\deg Q_m}{k}-1,\quad\phi=u|_\Sigma.\]
	Writing $\Delta_{\{p=0\}}+|A_{\{p=0\}}|^2$ in polar coordinate as $\partial_r^2+(n-2)r^{-1}\partial_r+r^{-2}(\Delta_\Sigma +|A_\Sigma|^2)$ on $\{p=0\}\cap\{\nabla p\neq0\}$, we have $\Delta_\Sigma \phi+|A_\Sigma|^2\phi=-\lambda \phi$ where $\lambda=\gamma^2-(n-3)\gamma$. For any $\Omega\subset\subset \Sigma$, we define
	\[\lambda_1(\Omega)=\inf_f \frac{\int_\Sigma |\nabla_\Sigma f|^2-|A_\Sigma|^2f^2}{\int_\Sigma f^2}\]
	where the $\inf$ is taken over all Lipschitz functions with compact support in $\Omega$. We also define 
	\begin{equation}\label{lambda1}
		\lambda_1(\Sigma)=\inf_{\Omega}\lambda_1(\Omega)
	\end{equation} where the $\inf$ is taken over all $\Omega\subset\subset \Sigma$ (see also \cite[(2.8)]{zhu2018first}). By the remark under \cite[Theorem 1.1]{zhu2018first}, $\Sigma$ satisfies the $\alpha$-structural hypothesis of \cite[Theorem 1.1]{zhu2018first} since it has a codimension 7 singular set by \cite{Fed70}. Thus by \cite[Theorem 1.1]{zhu2018first}, we have $\lambda_1(\Sigma)\le -(n-2)$.
	
	By Barta's theorem (see \cite{Barta} and also \cite[Lemma 1.36]{CM}) if there exists $\phi>0$ such that $(\Delta _\Sigma +|A_\Sigma|^2+\lambda)\phi=0$ on $\Sigma$ then for any $\Omega\subset\subset\Sigma$ we have $\lambda_1(\Omega)\ge \lambda$. Thus $\lambda_1(\Sigma)\ge \lambda$. Hence $\lambda\le -(n-2)$. We then obtain $\gamma^2-(n-3)\gamma+(n-2)\le0$. Thus 
	\[\mu_n^--1\le \gamma\le \mu_n^+-1.\] This shows that
	\[\mu_n^-\le\deg p+\frac{\deg Q_m}{k} \le \mu_n^+.\]
	
	We now show that equality on both sides cannot hold. Suppose we have equality on either side. Then we have $\lambda=-(n-2)$ and hence $\lambda_1(\Sigma)=-(n-2)$. By the rigidity part of \cite[Theorem 1.1]{zhu2018first} (see also \cite{perdomo2002first,wu1993new}), we have $\Sigma$ is a Clifford hypersurface. In particular, $\deg p=2$. This contradicts Theorem \ref{noisopara} that we are going to prove in the next section. There is no circular reasoning in between.
\end{proof}

\section{Further restrictions on the tangent cone of \texorpdfstring{$\graph P$}{graph P}}
\label{5}
As further consequences of Theorem \ref{struct}, we give some constraints on the tangent cone at infinity of the graph of a polynomial solution.

We first recall some concepts related to minimal cones.
\begin{definition}\label{sssm}
	A minimal hypercone $C\subset \R^n$ is called \textit{regular} if $\sing C\subset\{0\}$. Let $\Sigma=C\cap \SSS^{n-1}$ be the link of the cone and $\lambda_1$ be the first eigenvalue of $-(\Delta_\Sigma+|A_\Sigma|^2)$ as in (\ref{lambda1}). The minimal cone $C$ is called \textit{stable} (\textit{strictly stable}) if $(n-3)^2+4\lambda_1\ge 0$ ($>0$). There are two independent positive Jacobi fields on a stable cone $C$. If $C$ is strictly stable, they are $|x|^{-\gamma_\pm}\phi_1$. If $C$ is not strictly stable, they are $|x|^{-(n-3)/2}\phi_1$ and $|x|^{-(n-3)/2}(\log |x|)\phi_1$. Here $\phi_1$ is the first eigenfunction of $\Delta_\Sigma+|A_\Sigma|^2$ on $\Sigma$ and
	\begin{equation}\label{gamma}
		\gamma_{\pm}=\frac{n-3\pm\sqrt{(n-3)^2+4\lambda_1}}{2}.
	\end{equation}

	By \cite[Theorem 2.1]{simon1985area} if $C$ is a regular area-minimizing cone, then $\R^{n}\setminus C$ consists of two components $E_\pm$ and each $E_\pm$ is foliated uniquely by homotheties of a smooth minimal hypersurface $S_\pm\subset E_\pm$ called \textit{Hardt-Simon foliation}. There exist $R=R(C)>0,\alpha>0$ so that $S_\pm\setminus B_R(0)=\graph_C(v_\pm)$ and as $|x|\to\infty$
	\[ v_\pm(x)=w_\pm(x)+O(|x|^{-\gamma_{+}-\alpha})\] where $\pm w_\pm>0$ and $\Delta _Cw_\pm+|A_C|^2w_\pm=0$.
	
	By \cite[Theorem 3.1]{simon1985area}, a regular area-minimizing cone $C$ is called \textit{strictly minimizing} if $w_+$ and $w_-$ have slower decay near infinity.
\end{definition}

\begin{definition}
	An embedded hypersurface $\Sigma\subset \SSS^{l-1}$ is called \textit{isoparametric} if $\Sigma$ has constant principal curvature in $\SSS^{l-1}$. If  $\Sigma$ is isoparametric, so are its parallel hypersurfaces. We call $\Sigma$ and its parallel hypersurfaces an \textit{isoparametric family} of hypersurfaces in $\SSS^{l-1}$. For any isoparametric family, there is a unique hypersurface which is minimal. A regular minimal hypercone $C\subset\R^l$ is called \textit{isoparametric} if the link $C\cap \SSS^{l-1}$ is isoparametric and minimal. By \cite{Car,Mu1,Mu2}, hypersurfaces in an isoparametric family are the level sets of a polynomial $p$ in $\R^n$ restricted to $\SSS^{l-1}$. The polynomial $p$ is called a \textit{Cartan-M\"unzner polynomial} and satisfies the following equations
    $\Delta p=c|x|^{g-2}$ and $|\nabla p|^2=g^2|x|^{2g-2}$ where $g=\deg p$ and $c$ is some explicit constant depending on the family.
\end{definition}

 We first need the following information on the minimal hypersurface in an isoparametric family.
\begin{lemma}[\cite{Abr,Car,Mu1,Mu2,Wan}]\label{classisop}
	Let $p $ be a Cartan-M\"unzner polynomial on $\R^l$ and  $\Sigma$ be the minimal hypersurface in the isoparametric family associated to $p$. Then $\deg p$ is equal to the number of distinct principal curvatures of  $\Sigma\subset \SSS^{l-1}$. Moreover the following are all possibilities for the degree $\deg p$, the number of variables $l$ of $p$ and the norm squared of the second fundamental form $|A|^2$ of $\Sigma\subset \SSS^{l-1}$:
	
	(i) $\deg p =1$, $l\ge 2$, $|A|^2=0$,
	
	(ii) $\deg p =2$, $l\ge 4$, $|A|^2=l-2$,
	
	(iii) $\deg p =3$, $l=5,8,14,26$, $|A|^2=2(l-2)$,
	
	(iv) $\deg p =4$, $l\ge 6$ is even, $|A|^2=3(l-2)$,
	
	(v) $\deg p =6$, $l=8,14$, $|A|^2=5(l-2)$.
	
	Moreover, the cone $C$ over $\Sigma$ is area-minimizing if and only if $l\ge 4\deg p $ and $\Sigma$ is not $S^1\times S^5$ or $(SO(8)\times SO(2))/(SO(6)\times\Z_2)$. If $C$ is area-minimizing, then it is strictly stable and strictly minimizing.
\end{lemma}

\begin{proof}
	By \cite[Satz A and Satz 2]{Mu1}, $\deg p $ is equal to the number of distinct principal curvatures of the link and it can only be $1,2,3,4,6$. By \cite{Car}, for $\deg p \le 3$, the isoparametric hypersurfaces are homogeneous and $l$ can only take the values we mention. By \cite[Satz 1]{Mu1} if the $i$-th principal curvature has multiplicity $m_i$, then $m_i=m_{i+2}$ with index understood mod $\deg p $. Thus in case (iv), we have $l=2(m_1+m_2)+2\ge 6$ is even. For $\deg p =6$, by \cite{Abr}, $m_1=m_2$ can only be $1$ or $2$ and thus $l=3(m_1+m_2)+2=8$ or $14$. The values of $|A|^2$ follow from the exact values of principal curvatures of $\Sigma$  \cite[Satz 1]{Mu1}. The calculation can be found in \cite[Equations (48) (49) (50), and the formula below (53)]{TZ}. The area-minimizing and strictly minimizing properties are results of \cite{Wan}.
\end{proof}

\begin{remark}
    Given  a Cartan-M\"unzner polynomial $p $ on $\R^l$ and  $\Sigma$ the minimal hypersurface in the isoparametric family associated to $p$, we may not have $\Sigma=\{p=0\}\cap \SSS^{l-1}$. For example, $p(x)=(x_1^2+\cdots+x_r^2)-(x_{r+1}^2+\cdots+x_{r+s}^2)$ is a Cartan-M\"unzer polynomial on $\R^{l}$ for $l=r+s$. If $r\neq s$, then $\{p=0\}\cap \SSS^{l-1}$ has nonzero constant mean curvature in $\SSS^{l-1}$. However, there exists a homogeneous polynomial $q$ with the same degree as the Cartan-M\"unzer polynomial such that $\{q=0\}\cap \SSS^{l-1}=\Sigma$ as in the following lemma.
\end{remark}

\begin{lemma}\label{tilde.p}
    For any isoparametric minimal cone $C$, there exists a homogeneous polynomial $q$  such that $C=\{q=0\}$ and $\deg q=\deg p$ where $p$ is the Cartan-M\"unzner polynomial for the isoparametric family.
\end{lemma}
\begin{proof}
    Let $\Sigma\subset \SSS^{l-1}$ be the link of $C$. Let $p(\Sigma)=\{t\}$. If $\deg p=l$ is even, we define $q(x)=p(x)-t |x|^{l}$. Then $C=\{q=0\}$ and $\deg q=\deg p$.

    If $\deg p$ is odd, by Lemma \ref{classisop}, $\deg p=1$ or $3$.
    
    If $\deg p=1$, then $\{p=0\}$ is a hyperplane through the origin which is minimal.

If $\deg p = 3$, then $\{p = 0\}$ is minimal, which can be verified directly using the explicit expression given in \cite[Example~2]{TkClifford}. In both cases, we take $q=p$.
\end{proof}
\begin{theorem}\label{noisopara}
	Let $P=P_m+\cdots+P_1$ be a polynomial solution with $P_m=p ^kQ_m$ as in Theorem \ref{struct} (ii). Then $p$ cannot be a nonzero constant multiple of any polynomial given in Lemma \ref{tilde.p}.
	In particular we have $\deg p \ge 3$.
\end{theorem}

\begin{proof}
	Suppose for contradiction that $p $ is a nonzero constant multiple of a polynomial given in Lemma \ref{tilde.p} of degree $\ge 2$. Permuting the $x_i$'s if necessary, we assume that $x_1,\dots,x_l$ are all variables involved in $p$ for some $1\le l\le n$.
	
	We claim that $P_m=p ^kQ_m$ is independent of $x_{l+1},\dots,x_n$. Indeed, by Theorem \ref{struct} (i) and Lemma \ref{classisop}, the isoparametric minimal cone $\{p =0\}$ is strictly stable and strictly minimizing. By Lemma \ref{reg} and Theorem \ref{struct} (i), for $t\neq 0$, $\{P_m=t\}$ is a smooth minimal hypersurface lying on one side of the cone $\{p =0\}\times \R^{n-l}$. Moreover, the blow-down of $\{P_m=t\}$ is $\{p=0\}\times \R^{n-l}$ by Corollary \ref{convlscor}. By the Liouville theorem of \cite{ES}, $\{P_m=t\}$ is of the form $H\times \R^{n-l}$ where $H$ is a leaf of the Hardt-Simon foliation of $\{p =0\}$ in $\R^l$. Thus for all $t\in \R$, $\xi\in \R^{n-l}$, $\{P_m=t\}$ is invariant under $x\mapsto x+\xi$. This show that $P_m$ is independent of $x_{l+1},\dots,x_n$.
	
	Since $\{p =0\}$ is strictly minimizing, the Hardt-Simon foliation $\{P_m=t\}$ is a graph of $v$ over $\{p =0\}$ with $v$ decays like $|x|^{-\gamma_{-}}\phi_1$ near infinity. In particular, the velocity of the family of hypersurfaces $\{P_m=t^k\}$ is a multiple of $|x|^{-\gamma_{-}}\phi_1$. By Theorem \ref{struct} (iv), $\frac{1}{|\nabla (p Q_m^{1/k})|}$ is the velocity of the family $(\{P_m=t^{k}\})_{t\in (-\epsilon,\epsilon)}$ at $t=0$ near regular points of $\{P_m=0\}$. Then we have
	\begin{align}\label{Jacobi}
		\deg p -1+\frac{\deg Q_m}{k}=\frac{l-3-\sqrt{(l-3)^2+4\lambda_1}}{2}.
	\end{align}
	
	Equation (\ref{Jacobi}) implies
	\begin{equation}
		\label{crtn1} (l-3)^2+4\lambda_1\text{ is a perfect square integer.}
	\end{equation}
	To see this, since $|A|^2$ is an \textit{integer} for isoparametric minimal cones by Lemma \ref{classisop}, we have $\lambda_1=-|A|^2$ is also an integer. Since (\ref{Jacobi}) implies that the square root in (\ref{Jacobi}) is rational, it follows that the integer $(l-3)^2+4\lambda_1$ is a perfect square.
	
	By Theorem \ref{struct} (ii), we also have
	\begin{equation}
		\label{crtn2} k\text{ is odd and }\deg Q_m\text{ is even.}
	\end{equation}
	
	We claim that (\ref{Jacobi}), (\ref{crtn1}), (\ref{crtn2}) cannot hold simultaneously if $p $ is is a nonzero constant multiple of a polynomial given in Lemma \ref{tilde.p}. Indeed, we show this for each case in Lemma \ref{classisop}.
	
	\textbf{Case (i)}  $\deg p =1$. This is ruled out in Theorem \ref{struct} (ii).
	
	\textbf{Case (ii)} $\deg p =2$, $\lambda_1=-(l-2)$.
	
	In order for $\{p =0\}$ to be area-minimizing, we have $l\ge 8$ by Lemma \ref{classisop}. By (\ref{crtn1}), there is $a\in\N$ so that \[(l-3)^2-4(l-2)=a^2.\]
	This is equivalent to 
	\[(l-5-a)(l-5+a)=8.\]
	The only nonnegative integer solution to this equation is
	\[(l,a)=(8,1).\]
	By (\ref{Jacobi}), we have $k^{-1}\deg Q_m=1$. This contradicts (\ref{crtn2}).
	
	\textbf{Case (iii)} $\deg p =3$.
	
	In order for $\{p =0\}$ to be area-minimizing, we have $l=14$ or $l=26$ by Lemma \ref{classisop}.
	
	If $l=14$, then $\lambda_1=-24$ and thus $(l-3)^2+4\lambda_1=121-96=25$. By (\ref{Jacobi}) $k^{-1}\deg Q_m=1$. This contradicts (\ref{crtn2}).
	
	If $l=26$, then $\lambda_1=-48$ and $(l-3)^2+4\lambda_1=337$ is not a perfect square. This contradicts (\ref{crtn1}).
	
	\textbf{Case (iv)} $\deg p =4$. 
	
	In order for $\{p =0\}$ to be area-minimizing, we have $l\ge 16$ by Lemma \ref{classisop}.
	By (\ref{crtn1}), there is $a\in\N$ so that \[(l-3)^2-12(l-2)=a^2.\]
	This is equivalent to 
	\[(l-9-a)(l-9+a)=48.\]
	The nonnegative integer solutions to this equation are 
	\[(l,a)=(16,1),(17,4),(22,11).\]
	
	Since $l$ must be even by Lemma \ref{classisop}, the solution $(l,a)=(17,4)$ is not possible. 
	
	If $(l,a)=(16,1)$, from (\ref{Jacobi}) we have $k^{-1}\deg Q_m=3$. 
	
	If $(l,a)=(22,11)$, from (\ref{Jacobi}) we have $k^{-1}\deg Q_m=1$. Both violate (\ref{crtn2}).
	
	\textbf{Case (v)} $\deg p =6$.
	
	By Lemma \ref{classisop}, $\{p =0\}$ is not area-minimizing in this case. This finishes the proof of the main theorem.
	
	To show that $\deg p\ge 3$, we assume $\deg p=2$ such that $P_m=p^kQ_m$. By \cite[Observation 4]{Hsiang}, the minimal cone $C=\{p=0\}$ which is also the zero locus of a quadratic polynomial must be the Lawson's cone $C=C(S^r\times S^s)$. Thus $C$ is isoparametric and $C=\{q=0\}$ for $q$ given in Lemma \ref{tilde.p}. One can verify directly that $q$ is irreducible. Since $p,q$ are irreducible, we have $p=cq$ for $c\neq 0$ by Lemma \ref{null}. This contradicts what we just proved.
\end{proof}

We have the following corollary.
\begin{corollary}\label{noisopara.geometric}
	Let $P$ be a polynomial solution. Then the tangent cone at infinity of $\graph P$ cannot be $(C\times\R^{l-1})\times \R$ where $C$ is an isoparametric minimal cone and $l\ge 1$.
\end{corollary}
\begin{proof}
	Suppose the tangent cone at infinity of $\graph P$ is $(C\times\R^{l-1})\times \R$ where $C=\{q=0\}$ is an isoparametric minimal cone and $q$ is a polynomial on $\R^{n-l+1}$ given by Lemma \ref{tilde.p}. Viewing $q$ as a polynomial on $\R^n$, we can write $C\times\R^{l-1}=\{q=0\}$. By Theorem \ref{struct} (iii), $\{p=0\}=\{q=0\}$ where $p$ is the tangent cone factor of $P$. By Lemma \ref{null}, we have $p\mid q$. By Theorem \ref{noisopara}, $p\neq cq$ for any $c\neq0$ and $\deg p\ge 3$. Thus $\deg p<\deg q$ and $\deg q\ge 4$. By Lemma \ref{classisop} and \ref{tilde.p}, $\deg q=4$ or $6$. 
	
	If $\deg q=4$, then $\deg p=3$. Thus $q=pr$ where $r$ is a polynomial of degree 1. Then $\{q=0\}=\{p=0\}\cup\{r=0\}$ is the union of $\{p=0\}$ with a hyperplane. Since isoparametric minimal cones are regular, we must have $\{q=0\}$ is a hyperplane. This contradicts Lemma \ref{classisop} and \ref{tilde.p} which state that $\deg q=4$ is the number of distinct principal curvatures of $\{q=0\}\cap \SSS^{n-1}$.
	
	If $\deg q=6$ then by Lemma \ref{classisop} and \ref{tilde.p}, $C$ is not area-minimizing. This contradicts Theorem \ref{struct}.
\end{proof}

Our second restriction concerns polynomial solutions on $\R^8$ which is the lowest dimension that a polynomial solution could possibly exist.

\begin{theorem}\label{R8cubic}
	Let $P=P_m+\cdots+P_1$ be a polynomial solution on $\R^8$ with $P_m=p^kQ_m$ as in Theorem \ref{struct} (ii). Then $\deg p=3$ and $\{p=0\}\subset \R^8$ is a strictly stable, area-minimizing cone which is not strictly minimizing.
\end{theorem}

\begin{proof}
	We take $n=8$ throughout the proof. By Theorem \ref{n-2}
	\begin{equation}\label{rate}
		\deg p < \frac{n-1+\sqrt{(n-3)^2-4(n-2)}}{2}= 4.
	\end{equation}
	By Theorem \ref{noisopara}, we know that $\deg p \ge 3$. Thus $\deg p =3$. 
	
	Since $\{p=0\}$ is a $7$-dimensional area-minimizing cone, it is regular by \cite{Fed70}. We show that $\{p =0\}$ is strictly stable and not strictly minimizing. Indeed suppose $\{p =0\}$ were not strictly stable, then $\lambda_1=-\frac{25}{4}$. Since $u=\frac{1}{|\nabla (p Q_m^{1/k})|}$ is a positive Jacobi field on $\{p=0\}$ decaying without a $\log$ factor, we have $u=|x|^{-5/2}\phi_1$. In particular
	\[\deg p -1+k^{-1}\deg Q_m=\frac{5}{2}.\] Hence $k^{-1}\deg Q_m=\frac{1}{2}$. This contradicts the fact that $k$ is odd and $\deg Q_m$ is even by Theorem \ref{struct} (ii). Now suppose $\{p =0\}$ is strictly minimizing, then we can apply \cite[Theorem 5]{simon1989entire} to see that \[\deg P\le \gamma_{-}+1\le \frac{7}{2}\]where we used $\gamma_-\le \frac{5}{2}$ for $n=8$. This contradicts Theorem \ref{nodeg23} where cubic polynomial solutions are ruled out.
\end{proof}

The following corollary is a rephrasing of Theorem \ref{R8cubic}.
\begin{corollary}
	If $P$ is a polynomial on $\R^8$ solving the minimal surface equation and the tangent cone of $\graph P$ at infinity is $C\times\R$ where $C$ is strictly minimizing, then $P$ is an affine function.
\end{corollary}
\bmhead{Acknowledgements}

	The author would like to express his sincere gratitude to his advisor, Richard Schoen, for constant support and guidance. He would also like to thank Connor Mooney and Zhihan Wang for helpful discussions on this problem. Finally, he thanks the anonymous referees for many valuable suggestions.

\section*{Declarations}
\begin{itemize}
	\item Funding: The author did not receive support from any organization for the submitted work.
	\item Conflict of interest/Competing interests: The authors have no relevant financial or non-financial interests to disclose.
\end{itemize}



\begin{thebibliography}{9}
	\bibitem{Abr}
	Abresch, U.: Isoparametric hypersurfaces with four or six distinct principal curvatures. Math. Ann. \textbf{264}, 283-302 (1983)
	
	\bibitem{almgren1966some}
	Almgren, F.: Some interior regularity theorems for minimal surfaces and an extension of Bernstein's theorem. Ann. of Math. \textbf{84(1)}, 277-292 (1966)
	
	
	\bibitem{Barta}
	Barta, J.: Sur la vibration fondamentale d'une membrane. C. R. Acad. Sci. Paris \textbf{204}, 472-473 (1937)
	
	\bibitem{B}
	Bernstein, S.: \"Uber ein geometrisches Theorem und seine Anwendung auf die partiellen Differentialgleichungen vom elliptischen Typus. Math. Z. \textbf{26}, 551-558 (1927)
	
	\bibitem{bochnak2013real}
	Bochnak, J., Coste, M., Roy, M.: Real algebraic geometry. Springer Science \& Business Media \textbf{36} (2013)
	
	\bibitem{bombieri1969minimal}
	Bombieri, E., De Giorgi, E., Giusti, E.: Minimal cones and the Bernstein problem. Invent. Math. \textbf{7}, 243-268 (1969)
	
	\bibitem{bombieri1972harnack}
	Bombieri, E., Giusti, E.: Harnack's inequality for elliptic differential equations on minimal surfaces. Invent. Math. \textbf{15(1)}, 24-46 (1972)
	
	\bibitem{Car0}
	Cartan, E.: Familles de surfaces isoparam\'etriques dans les espaces \`a courbure constante. Ann. Mat. Pura Appl. \textbf{17}, 177-191 (1938)
	
	\bibitem{Car}
	Cartan, E.: Sur des familles remarquables d'hypersurfaces isoparam\'etriques dans les espaces sph\'eriques. Math. Z. \textbf{45}, 335-367 (1939)

    
	\bibitem{CM}
	Colding, T.H., Minicozzi II, W.P.: A Course in Minimal Surfaces. Amer. Math. Soc., Graduate Studies in Mathematics \textbf{121} (2011)
	
	\bibitem{de1965estensione}
	De Giorgi, E.: Una estensione del teorema di Bernstein. Ann. Scuola Norm. Sup. Pisa Cl. Sci. \textbf{19(1)}, 79-85 (1965)
	
	\bibitem{ES}
	Edelen, N., Sz\'ekelyhidi, G.: A Liouville-type theorem for cylindrical cones. Comm. Pure Appl. Math. \textbf{77}(8), 3557-3580 (2024).
	
	\bibitem{Fed}
	Federer, H.: Geometric Measure Theory. Classics in Mathematics, Springer, Berlin (2014)
	
	\bibitem{Fed70}
	Federer, H.: The singular sets of area minimizing rectifiable currents with codimension one and of area minimizing flat chains modulo two with arbitrary codimension. Ann. Math. \textbf{92}, 767-771 (1970)
	
	
	\bibitem{fleming1962oriented}
	Fleming, W.: On the oriented Plateau problem. Rend. Circ. Mat. Palermo \textbf{11(1)}, 69-90 (1962)
	
	\bibitem{GT}
	Gilbarg, D., Trudinger, N.S.: Elliptic Partial Differential Equations of Second Order. Grundlehren der Mathematischen Wissenschaften, vol. \textbf{224}. Springer, Berlin (1998)
	
	\bibitem{simon1985area}
	Hardt, R., Simon, L.: Area minimizing hypersurfaces with isolated singularities. J. Reine Angew. Math. \textbf{362}, 102-129 (1985)
	
	\bibitem{hartshorne1977algebraic}
	Hartshorne, R.: Algebraic Geometry. Graduate Texts in Mathematics, vol. 52. Springer-Verlag, New York (1977)
	
	\bibitem{Hsiang}
	Hsiang, W.-Y.: Remarks on closed minimal submanifolds in the standard Riemannian $m$-sphere. J. Differential Geom. \textbf{1}, 257-267 (1967)
	
	\bibitem{lin_1987} 
	Lin, F.-H.: Minimality and stability of minimal hypersurfaces in $\mathbb{R}^n$. Bull. Aust. Math. Soc. \textbf{2}, 209-214 (1987)
	
	\bibitem{morgan2016gmt}
	Morgan, F.: Geometric Measure Theory: A Beginner's Guide, Academic Press, Amsterdam (2016)
	
	\bibitem{Mu1}
	M\"unzner, H. F.: Isoparametric Hyperfl\"achen in Sph\"aren, I. Math. Ann. \textbf{251}, 57-71 (1980)
	
	\bibitem{Mu2}
	M\"unzner, H. F.: Isoparametric Hyperfl\"achen in Sph\"aren, II. Math. Ann. \textbf{256}, 215-232 (1981)
	
	
	\bibitem{NTV}
	Nadirashvili, N., Tkachev, V., Vladut, S.: Nonlinear elliptic equations and nonassociative algebras. Vol. 200. American Mathematical Soc. (2014)
	
	\bibitem{perdomo2002first}
	Perdomo, O.: First stability eigenvalue characterization of Clifford hypersurfaces. Proc. Amer. Math. Soc. \textbf{130(11)}, 3379-3384 (2002)
	
	\bibitem{Sim}
	Simon, L.: Lectures on geometric measure theory. The Australian National University, Mathematical Sciences Institute, Centre for Mathematics \& its Applications (1983)
	
	\bibitem{simon1983survey}
	Simon, L.: Survey lectures on minimal submanifolds. Ann. Math. Stud. \textbf{103}, 3-52 (1983)
	
	\bibitem{simon1987asymptotic}
	Simon, L.: Asymptotic behaviour of minimal graphs over exterior domains. Ann. Henri Poincar\'e \textbf{4(3)}, 231-242 (1987)
	
	\bibitem{simon1989entire}
	Simon, L.: Entire solutions of the minimal surface equation. J. Differ. Geom. \textbf{30(3)}, 643-688 (1989)
	
	\bibitem{simon1997minimal}
	Simon, L.: The minimal surface equation. In: Geometry V, pp. 239-266 (1997)
	
	\bibitem{simons1968minimal}
	Simons, J.: Minimal varieties in Riemannian manifolds. Ann. of Math. \textbf{88}, 62-105 (1968)
	
	\bibitem{TZ}
	Tang, Z., Zhang, Y.: Minimizing cones associated with isoparametric foliations. J. Differ. Geom. \textbf{115(2)}, 367-393 (2020)
	
	\bibitem{TkClifford}
	Tkachev, V.: Minimal cubic cones via Clifford algebras. Complex Anal. Oper. Theory \textbf{4}, 685-700 (2010)
	
	\bibitem{tkachev2010classification}
	Tkachev, V.: On a classification of minimal cubic cones in $\mathbb{R}^{n}$. Preprint (2010). Available at: \url{https://arxiv.org/abs/1009.5409}
	
	\bibitem{tkachev2016non}
	Tkachev, V.: On the non-vanishing property for real analytic solutions of the $p$-Laplace equation. Proc. Amer. Math. Soc. \textbf{144(6)}, 2375-2382 (2016)
	
	\bibitem{Wan}
	Wang, Q.-M.: On a class of minimal hypersurfaces in $\mathbb{R}^n$. Math. Ann. \textbf{298}, 207-251 (1994)
	
	\bibitem{wu1993new}
	Wu, C.: New characterizations of the Clifford tori and the Veronese surface. Arch. Math. (Basel) \textbf{61(3)}, 277-284 (1993)
	
	\bibitem{Yau}
	Yau, S.-T. (ed.): Seminar on differential geometry. No. 102. Princeton University Press (1982)
	
	\bibitem{zhu2018first}
	Zhu, J.: First stability eigenvalue of singular minimal hypersurfaces in spheres. Calc. Var. Partial Differ. Equ. \textbf{57(5)}, 1-13 (2018)
	
\end{thebibliography}
\end{document}